\documentclass[a4paper,11pt]{article}

\usepackage{amsmath}
\usepackage{amsthm}
\usepackage{amsfonts}
\usepackage{amssymb}
\usepackage{amscd}
\usepackage[pdftex]{graphicx,xcolor,hyperref}
\usepackage{tikz}
\usetikzlibrary{intersections,calc,arrows.meta,cd}
\usepackage{mathrsfs}
\usepackage{here}

\newtheorem{thm}{Theorem}[section]

\newtheorem{defn}[thm]{Definition}
\newtheorem{lem}[thm]{Lemma}
\newtheorem{prop}[thm]{Proposition}
\newtheorem{rem}[thm]{Remark}
\newtheorem{exam}[thm]{Example}

\numberwithin{equation}{section}

\begin{document}

\title{\bf \LARGE  On \textit{q}-deformed Farey sum and a homological interpretation of \textit{q}-deformed real quadratic irrational numbers}
\author{Xin Ren}
\date{}

\maketitle

\begin{abstract}
The left and right \textit {q}-deformed rational numbers were introduced by Bapat, Becker and Licata via regular continued fractions, and they gave a homological interpretation for left and right \textit {q}-deformed rational numbers. In the present paper, we focus on negative continued fractions and defined left \textit {q}-deformed negative continued fractions. We give a formula for computing the \textit {q}-deformed Farey sum of the left \textit {q}-deformed rational numbers based on it. We use this formula to give a combinatorial proof of the relationship between the left \textit {q}-deformed rational number and the Jones polynomial of the corresponding rational knot which was proved by Bapat, Becker and Licata using a homological technique. Finally, we combine their work and the \textit {q}-deformed Farey sum, and give a homological interpretation of the \textit {q}-deformed Farey sum. We also give an approach to finding a relationship between real quadratic irrational numbers and homological algebra.
\end{abstract}
keywords: continued fractions, \textit {q}-deformed Farey sum, Jones polynomial, real quadratic irrational numbers, 2-Calabi–Yau category
\footnote[0]{MSC Classification(2020): 11A55, 05A30, 57K14, 18G80}


\section{Introduction}

\noindent
The notion of $q$-deformed rational numbers \cite{MO1} was introduced by Morier-Genoud and Ovsienko based on some combinatorial properties of rational numbers. They further extended this notion to arbitrary real numbers \cite{MO2} by some number-theoretic properties of irrational numbers. These works are related to many directions including Jones polynomial of rational knots \cite{KW, KR, NT, MO1}, Teichm\"uller spaces \cite{VL}, the Markov-Hurwitz approximation theory \cite{K, SM, LMOV, XR}, the modular group and the Picard group \cite{LM, O}, combinators of posets \cite{TBEC1, OEK1, OEK2} and triangulated category~\cite{BBL}.

\indent
For a formal parameter $q$ and an irreducible fraction $\displaystyle\frac{r}{s}$, as an enhancement of $q$-deformed rational numbers, Bapat, Becker and Licata defined left $q$-deformed rational number $\displaystyle\left[\frac{r}{s}\right]^{\flat}_q$ and right $q$-deformed rational number $\displaystyle\left[\frac{r}{s}\right]^{\sharp}_q$ via the regular continued fractions of $\displaystyle\frac{r}{s}$, and the right $q$-deformed rational number $\displaystyle\left[\frac{r}{s}\right]^{\sharp}_q$ is exactly $q$-deformed rational number $\displaystyle\left[\frac{r}{s}\right]_q$ considered by Morier-Genoud and Ovsienko.  Following \cite{MO1} and \cite{BBL}, the right $q$-deformed rational numbers can be expressed by the right $q$-deformed regular or negative continued fraction expansions, and the left $q$-deformed rational numbers can be expressed by the left $q$-deformed regular continued fraction expansions. These $q$-deformations of the fractions are rational expressions in the variable $q$ with integer coefficients. Such as \cite{MO1, LM}, and so on, it may be more convenient from the perspective of the negative continued fraction expansion when we consider some properties of left and right $q$-deformed rational numbers and their applications. In particular, the formula for the right $q$-deformed Farey sum based on the negative continued fraction is more concise \cite[Section 2]{MO1}. This induces us to consider the $q$-deformed Farey sum of the left $q$-deformed rational numbers. In the present paper, we define the left $q$-deformed negative continued fraction expansion. Then we give a formula for computing the $q$-deformed Farey sum of the left $q$-deformed rational numbers based on negative continued fraction (see Theorem \ref{THM-F-S-L}). 

\indent
As an application of the right $q$-deformed rational numbers, given a rational number $\displaystyle\frac{r}{s}$, we can use the numerator and denominator of $\displaystyle\left[\frac{r}{s}\right]^{\sharp}_q$ to represent the Jones polynomial of the rational knot to which $\displaystyle\frac{r}{s}$ corresponds \cite[Proposition A.1]{MO1}. On the other hand, Bapat, Becker and Licata prove that the Jones polynomial for the rational knot corresponding to $\displaystyle\frac{r}{s}$ can be represented by just the numerator of $\displaystyle\left[\frac{r}{s}\right]^{\flat}_q$ \cite[Theorem A.3]{BBL} by considering a homological interpretation of $\displaystyle\left[\frac{r}{s}\right]^{\flat}_q$ and $\displaystyle\left[\frac{r}{s}\right]^{\sharp}_q$. Considering the zigzag algebra on the $A_2$ quiver, we can obtain a triangulated category $\mathcal{C}_2$ called $2$-Calabi–Yau category associated to the $A_2$ quiver \cite[Section 4]{BBL}. For spherical objects on $\mathcal{C}_2$, Bapat, Becker and Licata defined two functions, denoted as $\mathrm{occ}_q$ and $\overline{\mathrm{hom}}_q$, and they proved that $\displaystyle\left[\frac{r}{s}\right]^{\flat}_q$ and $\displaystyle\left[\frac{r}{s}\right]^{\sharp}_q$ can be expressed in terms of $\mathrm{occ}_q$ and $\overline{\mathrm{hom}}_q$, respectively \cite[Theorems 3.7 and 3.8]{BBL}. There are two questions worth considering. Can we give a combinatorial proof of \cite[Theorem A.2]{BBL} without homology techniques? Can we give a homological interpretation of the $q$-deformed irrational numbers defined in \cite{MO2}? In the present paper, we apply Theorem \ref{THM-F-S-L} to give a combinatorial proof of \cite[Theorem A.3]{BBL} without using homology techniques (see Theorem \ref{J-L-Q-R}). Then, we combine the homological interpretation of the left and right $q$-deformed rational numbers and the $q$-deformed Farey sum, and give a homological interpretation of the $q$-deformed Farey sum (see Propositions \ref{CTFS-RandL} and \ref{CTFS-OB}). We also apply the results in \cite[Theorems 4.7 and 4.8]{BBL} to real quadratic irrational numbers with periodic type (see Theorem~\ref{RQIN-occ_q_and_hom_q}). 

\indent 
This paper is organized into the following sections. In Section \ref{Sec2}, we first recall some definitions related to the left and right $q$-deformed rational numbers, including the (right) \textit{q-deformed Euler continuants}, which were introduced by Morier-Genoud and Ovsienko \cite{MO1}. Similarly, we define the left $q$-deformed negative continued fractions and left $q$-deformed Euler continuants. We prove that the left $q$-deformed negative continued fractions and the left $q$-deformed regular continued fractions are equal.
  In Section \ref{Sec3}, we give $q$-deformed Farey sum of left $q$-rational numbers based on negative continued fraction, and derive a weighted triangulation and $q$-deformed Farey tessellation corresponding to left $q$-deformed rational numbers. In Section \ref{Sec4}, we give a new proof of \cite[Theorem A.2]{BBL} as an application of $q$-deformed Farey sum of left $q$-deformed rational numbers by induction through the length of the negative continued fraction expansion. 
In Section \ref{Sec5}, we first combine the results of \cite[Theorems 4.7 and 4.8]{BBL} with the $q$-deformed Farey sum to give a homological interpretation of the $q$-deformed Farey sum. Then we consider a real quadratic irrational number with periodic type, its $q$-deformation can also be expressed in a special form that is related to the $q$-deformation rational number that approximates it \cite{LM}. Based on the results of these $q$-deformations, we give the relations between real quadratic irrational numbers and homological algebra.


\section{$q$-deformed continued fractions and $q$-deformed Euler continuants}\label{Sec2}
\noindent
In this section, we first briefly review some definitions related to left and right $q$-deformed rational numbers (see \cite{BBL} and \cite{MO1} for details). We define the left $q$-deformed negative continued fraction expansion and introduce the left $q$-deformed Eulerian continuants. We simply check that the left $q$-deformed negative continued fraction expansion is indeed consistent with the left $q$-deformed regular continued fraction expansion.
 
\subsection{Left and right $q$-integers and $q$-deformed rational numbers}
It is well-known that an irreducible fraction $\displaystyle\frac{r}{s} \in \mathbb{Q}\cup \{\infty\}$ has unique regular and negative continued fraction expansions as follows:\\
\[ \displaystyle
\frac{r}{s}=a_1+\frac{1}{a_2+\cfrac{1}{\ddots+\cfrac{1}{a_{2m}}}}=c_1-\frac{1}{c_2-\cfrac{1}{\ddots-\cfrac{1}{c_k}}}
\]
\\
with $a_1 \in \mathbb{Z}$, $a_i \in \mathbb{Z}\setminus\{0\} $ ($i\geq2$), and $c_1 \in \mathbb{Z} \setminus\{0\} $ ($i\geq2$) and $c_j \in \mathbb{Z}\setminus\{-1,0,1\} $ ($j\geq2$). If $\displaystyle\frac{r}{s}$ is negative, then $a_1,\ldots, a_{2m}$ and $c_1,\ldots, c_{k}$ are negative, and if $\displaystyle\frac{r}{s}$ is positive, then $a_1,\ldots, a_{2m}$ and $c_1,\ldots, c_{k}$ are positive. We denote this expansion by $[a_1,\ldots,a_{2m}]$ and $[[c_1,\ldots, c_{k}]]$, respectively. As special cases, the regular and negative continued fraction expansions of $0$ and $\infty$  ($\displaystyle\infty:=\frac{1}{0}$)  are $[-1,1]$, $[[1,1]]$ and empty expansion $[\; ]$, $[[\; ]]$, respectively.

\medskip

\indent
We consider the following three matrices.
$\displaystyle
\sigma_1:=\begin{pmatrix}
1 & -1 \\
0 & 1 \\
\end{pmatrix}, \ \ \ \ 
\sigma_2:=\begin{pmatrix}
1 & 0 \\
1 & 1 \\
\end{pmatrix}, \ \ \ \
S:=~\begin{pmatrix}
0 & -1 \\
1 & 0 \\
\end{pmatrix}.
$
We know that the modular group $\mathrm{PSL}_2(\mathbb{Z})$ can be generated
by $\{\sigma_1,\sigma_2\}$ or $\{\sigma_1,S\}$. The modular group $\mathrm{PSL}_2(\mathbb{Z})$ acts on $\displaystyle \mathbb{Q}\cup \{\infty\}$ by the fractional linear transformation:
\[
\displaystyle
\begin{pmatrix}
a & b \\
c & d \\
\end{pmatrix}(x)=\frac{ax+b}{cx+b},\]
where
$\displaystyle
\begin{pmatrix}
a & b \\
c & d \\
\end{pmatrix} \in \mathrm{PSL}_2(\mathbb{Z})$, $\displaystyle x\in \mathbb{Q}\cup \{\infty\}$. Then a rational number $\displaystyle \frac{r}{s}=[a_1,\ldots,a_{2m}]=[[c_1,\ldots,c_k]]$ can be expressed by the following formulas:
\begin{equation}\label{R_r_over_s}
\displaystyle
\frac{r}{s}=\sigma_1^{-a_1}\sigma_2^{a_2}\sigma_1^{-a_3}\sigma_2^{a_4}\cdots \sigma_1^{-a_{2m-1}} \sigma_2^{a_{2m}}(\infty),
\end{equation}
\begin{equation}
\displaystyle
\frac{r}{s}=\sigma_1^{-c_1}S\sigma_1^{-c_2}S\cdots \sigma_1^{-c_k}S(\infty).
\end{equation}

\begin{defn}[\cite{BBL}]
{\normalfont
Let $q$ be a formal parameter. For a rational number $\displaystyle\frac{r}{s}=[a_1,\ldots,a_{2m}]$, 
we denote by $\mathrm{PSL}_{2,q}(\mathbb{Z})$ the subgroup of group $\mathrm{GL}_{2}(\mathbb{Z}[q^{\pm 1}])$ generated by the following two elements:
\[\displaystyle
\sigma_{1,q}=\begin{pmatrix}
q^{-1} & -q^{-1} \\
0 & 1 \\
\end{pmatrix}, \ \ \ \ 
\sigma_{2,q}=\begin{pmatrix}
1 & 0 \\
1 & q^{-1} \\
\end{pmatrix}. \ \ \ \
\]
Then the \textit {right} $q$-deformed rational number is
\[\displaystyle
\left[\displaystyle\frac{r}{s}\right]^{\sharp}_q=\sigma_{1,q}^{-a_1}\sigma_{2,q}^{a_2}\sigma_{1,q}^{-a_3}\sigma_{2,q}^{a_4}\cdots \sigma_{1,q}^{-a_{2m-1}} \sigma_{1,q}^{a_{2m}}\left(\infty\right),
\]
and the \textit {left} $q$-deformed rational number is

\[\displaystyle
\left[\displaystyle\frac{r}{s}\right]_q^{\flat}=\sigma_{1,q}^{-a_1}\sigma_{2,q}^{a_2}\sigma_{1,q}^{-a_3}\sigma_{2,q}^{a_4}\cdots \sigma_{1,q}^{-a_{2m-1}} \sigma_{1,q}^{a_{2m}}\left(\frac{1}{1-q}\right).
\]
}
\end{defn}

\subsection{Left $q$-deformed negative continued fractions}

\begin{defn}[\cite{BBL}]
{\normalfont
Let $q$ be a formal parameter.  We consider an integer $n$, the following two rational forms $[n]^{\flat}_q$ and $[n]^{\sharp}_q$ in $q$ are called the right $q$-integer of $n$ and the left $q$-integer of $n$, respectively.
\[
\displaystyle
[n]^{\sharp}_q:=\frac{1-q^n}{1-q}, \quad
\displaystyle
[n]^{\flat}_q:=\frac{1-q^{n-1}+q^n-q^{n+1}}{1-q}.
\]
}
\end{defn}

\begin{rem}\label{REM-flat-n}
{\normalfont
Suppose that $m, n \in \mathbb{Z}$. It can be easy to check that the right $q$-integers and left $q$-integers satisfy the following properties.
\begin{center}
\begin{enumerate}
\item[(i)] $\left[n\right]^{\sharp}_q=\left[n\right]^{\flat}_q+q^{n-1}-q^n$;
\item[(ii)] $\left[m+n\right]^{\sharp}_q=\left[m\right]^{\sharp}_q+q^m[n]^{\sharp}_q=\left[n\right]^{\sharp}_q+q^n\left[m\right]^{\sharp}_q=\left[n+n\right]^{\sharp}_q$, \\ $\left[m+n\right]^{\flat}_q=\left[m\right]^{\sharp}_q+q^m\left[n\right]^{\flat}_q=\left[n\right]^{\sharp}_q+q^n\left[m\right]^{\flat}_q=\left[n+m\right]^{\flat}_q$;
\item[(iii)] $\left[-n\right]^{\sharp}_q=-q^{-1}\left[n\right]^{\sharp}_{q^{-1}}$, $\left[-n\right]^{\flat}_q=-q^{-1}\left[n\right]^{\flat}_{q^{-1}}$;
\item[(iv)]  $\displaystyle q^n\left[n\right]^{\sharp}_{q^{-1}}=q\left[n\right]^{\sharp}_q$, $\displaystyle q^n(\left[n\right]^{\flat}_{q^{-1}}-\left[0\right]^{\flat}_{q^{-1}})=q(\left[n\right]^{\flat}_q-\left[0\right]^{\flat}_q)$.
\end{enumerate}
\end{center}
}
\end{rem}

Suppose that $\displaystyle\frac{r}{s}=[a_1,\ldots,a_{2m}]=[[c_1,\ldots,c_k]]$. From \cite{MO1} and \cite{BBL}, the right $q$-deformed rational number $\displaystyle\left[\frac{r}{s}\right]^{\sharp}_q$ has both the following $q$-deformed positive and negative continued fraction expansions.

\begin{equation}\label{eqn-1.3}
\displaystyle\left[\frac{r}{s}\right]^{\sharp}_q=[a_1, a_2, \ldots , a_{2m}]^{\sharp}_q :=
[a_1]^{\sharp}_q  
+\cfrac{ q^{a_1}   }{ [a_2]^{\sharp}_{ q^{-1} }  + \cfrac{ q^{-a_2}  }
{  [a_3]^{\sharp}_q  + \cfrac{q^{a_3} }
{ [a_4]^{\sharp}_{q^{-1} } + \cfrac{ q^{-a_4} }{  \cfrac{\ddots }
{ [a_{2m-1}]^{\sharp}_q+ \cfrac{ q^{a_{2m-1}} }{ [a_{2m}]^{\sharp}_{q^{-1}}  }}}}}}, 
\end{equation}

\noindent 

\begin{equation}\label{eqn-1.4}
\displaystyle\left[\frac{r}{s}\right]^{\sharp}_q=[[c_1, c_2, \ldots , c_{k}]]^{\sharp}_q :=
[c_1]^{\sharp}_q  
-\cfrac{ q^{c_1 -1}   }{ [c_2]^{\sharp}_{ q }  - \cfrac{ q^{c_2 -1}  }
{  [c_3]^{\sharp}_q  - \cfrac{q^{c_3-1 }}
{ [c_4]^{\sharp}_{q } - \cfrac{ q^{c_4 -1} }{  \cfrac{\ddots }
{ [c_{k-1}]^{\sharp}_q- \cfrac{ q^{c_{k-1}-1} }{ [c_{k}]^{\sharp}_{q}  }}}}}} . 
\end{equation}

For the left $q$-deformed rational number $\displaystyle\left[\frac{r}{s}\right]^{\flat}_q$, Bapat,Becker and Licata proved that the right $q$-deformed rational number $\displaystyle\left[\frac{r}{s}\right]^{\flat}_q$ has a $q$-deformed positive continued fraction expansion \cite{BBL} as follows: 

\begin{equation}\label{eqn-1.5}
\displaystyle\left[\frac{r}{s}\right]^{\flat}_q=[a_1, a_2, \ldots , a_{2m}]^{\flat}_q :=
[a_1]^{\sharp}_q  
+\cfrac{ q^{a_1}   }{ [a_2]^{\sharp}_{ q^{-1} }  + \cfrac{ q^{-a_2}  }
{  [a_3]^{\sharp}_q  + \cfrac{q^{a_3} }
{ [a_4]^{\sharp}_{q^{-1} } + \cfrac{ q^{-a_4} }{  \cfrac{\ddots }
{ [a_{2m-1}]^{\sharp}_q+ \cfrac{ q^{a_{2m-1}} }{ [a_{2m}]^{\flat}_{q^{-1}}  }}}}}}.
\end{equation}

\medskip

\noindent
Similarly to the formula \eqref{eqn-1.4}, we can also define the left $q$-deformed negative continued fraction expansion as follows:

\begin{defn}[left $q$-deformation of negative continued fraction expansion]
{\normalfont
\begin{equation}\label{eqn-1.7}
[[c_1, c_2, \ldots , c_{k}]]^{\flat}_q :=
[c_1]^{\sharp}_q  
-\cfrac{ q^{c_1 -1}}{ [c_2]^{\sharp}_{ q }  - \cfrac{ q^{c_2 -1}  }
{  [c_3]^{\sharp}_q  - \cfrac{q^{c_3-1 }}
{ [c_4]^{\sharp}_{q } - \cfrac{ q^{c_4 -1} }{  \cfrac{\ddots }
{ [c_{k-1}]^{\sharp}_q- \cfrac{ q^{c_{k-1}-1} }{ [c_{k}]^{\flat}_{q}  }}}}}} . 
\end{equation}
Note that it differs from the right $q$-deformed negative continuous fraction expansion only in the last term.
}
\end{defn}

As in the case of right $q$-deformation, we have the following conclusion for the case of left $q$-deformation.

\begin{thm}\label{Thm1}
If a rational number $\displaystyle\frac{r}{s}$ is given in the form 
\noindent 
$\displaystyle\frac{r}{s}=[a_{1}, \ldots , a_{2m}]=[[c_{1}, \ldots , c_{k}]]$, then 

\medskip

\noindent 
\begin{equation}\label{eqn-1.8}
[a_{1}, \dots , a_{2m}]^{\flat}_{q}=[[c_{1}, \ldots , c_{k}]]^{\flat}_{q}.
\end{equation}
\end{thm}
We will prove this formula in Section \ref{POThm1}.

\medskip

\noindent
By \cite{MO1} and \cite{BBL}, the left and right $q$-rationals can be expressed by the following quotients of two polynomial in $q$ with integer coefficients as follow:
\[\displaystyle
\left[\frac{r}{s}\right]^{\sharp}_q=\frac{\mathcal{R}^{\sharp}(q)}{\mathcal{S}^{\sharp}(q)}, \ \ \ \
\left[\frac{r}{s}\right]^{\flat}_q=\frac{\mathcal{R}^{\flat}(q)}{\mathcal{S}^{\flat}(q)},
\]where $\mathcal{R}^{\dag}(q)$ and $\mathcal{S}^{\dag}(q)$ are coprime, and $\mathcal{R}^{\dag}(q), \ \mathcal{S}^{\dag}(q) \in \mathbb{Z}[q]$ are monic polynomials ($\dag \in \{ \sharp, \flat \}$).
In particular, we have
\[
\displaystyle
\left[\frac{0}{1}\right]^{\sharp}_q=\frac{0}{1}, \ \ \ 
\left[\frac{0}{1}\right]^{\flat}_q=\frac{1-q^{-1}}{1}; \ \ \ 
\left[\infty\right]^{\sharp}_q=\frac{1}{0}, \ \ \ 
\left[\infty\right]^{\flat}_q=\frac{1}{1-q}. \ \ \ 
\]

\medskip


\subsection{The left and right $q$-deformed Euler continuants}\label{TQEC}
\begin{defn}[right \textit{q}-deformed Euler continuants]
{\normalfont
\begin{equation}
\displaystyle
E^{\sharp}_k(c_1,\ldots,c_k)_q:=
\begin{vmatrix}
[c_1]^{\sharp}_q & q^{c_1-1} & & & \\
1 & [c_2]^{\sharp}_q & q^{c_2-1} & & \\
 & \ddots & \ddots & \ddots &  \\
 &        &   1    &   [c_{k-1}]^{\sharp}_q     & q^{c_{k-1}-1} \\
 &        &        &       1           & [c_k]^{\sharp}_q  \\
\end{vmatrix}
\end{equation}
where $c_i$'s are integers, and for convenience, we set $E_0^{\sharp}()=1$ and $E_{-1}^{\sharp}()=0$.
}
\end{defn}

For the numerators and denominators of the right $q$-deformed rational numbers, we have the following conclusion.

\begin{prop}[{\cite[Proposition 5.3]{MO1}}]\label{N_and_D_and_q-Euler}
For a rational number $\displaystyle\frac{r}{s}=[[c_1,\ldots,c_k]]$, we have
\[
\displaystyle
\mathcal{R}^{\sharp}(q)= E^{\sharp}_k(c_1,\ldots,c_k)_q, \ \ \ \ \ \ 
\mathcal{S}^{\sharp}(q)= E^{\sharp}_{k-1}(c_2,\ldots,c_k)_q.
\]
\end{prop}

\begin{defn}[left \textit{q}-deformed Euler continuants]
{\normalfont
\begin{equation}\label{E-flat}
\displaystyle
E^{\flat}_k(c_1,\ldots,c_k)_q:=
\begin{vmatrix}
[c_1]^{\sharp}_q & q^{c_1-1} & & & \\
1 & [c_2]^{\sharp}_q & q^{c_2-1} & & \\
 & \ddots & \ddots & \ddots &  \\
 &        &   1    &   [c_{k-1}]^{\sharp}_q     & q^{c_{k-1}-1} \\
 &        &        &       1           & [c_k]^{\flat}_q  \\
\end{vmatrix}
\end{equation}
where $c_i$'s are integers, and for convenience, we set $E_0^{\flat}()=1$ and $E_{-1}^{\flat}()=1-q$.
}
\end{defn}

\medskip

By the definition of $E^{\sharp}_{k}(c_1,\dots,c_k)_q$, we know that for $c_l>2$, we have
\begin{equation}
\displaystyle
q(E^{\sharp}_{l+1}(c_1,\dots,c_l,2)_q-E^{\sharp}_{l}(c_1,\ldots,c_l)_q)=E^{\sharp}_{l+2}(c_1,\dots,c_l,2,2)_q-E^{\sharp}_{l+1}(c_1,\dots,c_l,2)_q.\notag
\end{equation}
For a non-negative integer $h$, by induction, it can be checked that 
\begin{equation}\label{E-flat-2}
\displaystyle
q^h(E^{\sharp}_{l+1}(c_1,\dots,c_l,2)_q-E^{\sharp}_{l}(c_1,\ldots,c_l)_q)=E^{\sharp}_{l+h+1}(c_1,\dots,c_l,2^{(h+1)})_q-E^{\sharp}_{l+h}(c_1,\ldots,c_l,2^{(h)})_q.
\end{equation}
Moreover, by the definition of $E^{\sharp}_{k}(c_1,\dots,c_k)_q$ and $E^{\flat}_{k}(c_1,\dots,c_k)_q$, we can know that
\begin{align}\label{E-flatE}
E^{\flat}_{k}(c_1,\dots,c_k)_q&=[c_k]^{\flat}_qE^{\sharp}_{k-1}(c_1,\dots,c_{k-1})_q-q^{c_{k-1}-1}E^{\sharp}_{k-2}(c_1,\dots,c_{k-2})_q\notag \\
&=E^{\sharp}_{k}(c_1,\dots,c_k)_q-q^{c_{k}-1}(1-q)E^{\sharp}_{k-1}(c_1,\dots,c_{k-1})_q.
\end{align}

\medskip


\subsection{Proof of Theorem \ref{Thm1}}\label{POThm1}

Before we prove Theorem \ref{Thm1}, let us prove the following proposition.
\begin{prop}
Consider the element $\displaystyle S_q:=\begin{pmatrix}
0 & -q^{-1} \\
1 & 0 \\
\end{pmatrix}$ in $\mathrm{PSL}_{2,q}(\mathbb{Z})$, then 
\begin{equation}
\displaystyle
\sigma_{1,q}^{-c_1}S_q\sigma_{1,q}^{-c_2}S_q\cdots \sigma_{1,q}^{-c_k}S_q\left(\displaystyle\frac{1}{1-q}\right)=[[c_1,\ldots,c_k]]^{\flat}_q.
\end{equation}
\end{prop}

\noindent
\begin{proof}
By Proposition 4.3 of \cite{MO1} and Proposition \ref{N_and_D_and_q-Euler}, one has 
\[
\displaystyle
\sigma_{1,q}^{-c_1}S_q\sigma_{1,q}^{-c_2}S_q\cdots \sigma_{1,q}^{-c_k}S_q=\begin{pmatrix}
E^{\sharp}_k(c_1,\ldots,c_k)_q & -q^{c_k-1}E^{\sharp}_{k-1}(c_1,\ldots,c_{k-1})_q \\
E^{\sharp}_{k-1}(c_2,\ldots,c_k)_q & -q^{c_k-1}E^{\sharp}_{k-2}(c_2,\ldots,c_{k-1})_q\\
\end{pmatrix}.
\]
We view the left $q$-rational $\displaystyle\left[\infty\right]^{\flat}_q=\frac{1}{1-q}$ as a vector in the projective space. Note that
\[
\displaystyle
E^{\sharp}_k(c_1,\ldots,c_k)_q=[c_k]^{\sharp}_qE^{\sharp}_{k-1}(c_1,\ldots,c_{k-1})_q-q^{c_{k-1}-1}E^{\sharp}_{k-2}(c_1,\ldots,c_{k-2})_q,
\]
and by Remark \ref{REM-flat-n}, we have
\\
\begin{align}
\displaystyle
&\sigma_{1,q}^{-c_1}S_q\sigma_{1,q}^{-c_2}S_q\cdots \sigma_{1,q}^{-c_k}S_q\begin{pmatrix}
1\\
1-q\\
\end{pmatrix}\notag\\ =\ &\begin{pmatrix}
E^{\sharp}_k(c_1,\ldots,c_k)_q-q^{c_k-1}E^{\sharp}_{k-1}(c_1,\ldots,c_{k-1})_q+q^{c_k}E^{\sharp}_{k-1}(c_1,\ldots,c_{k-1})_q\\
E^{\sharp}_{k-1}(c_2,\ldots,c_k)_q-q^{c_k-1}E^{\sharp}_{k-2}(c_2,\ldots,c_{k-1})_q+q^{c_k}E^{\sharp}_{k-2}(c_2,\ldots,c_{k-1})_q\\
\end{pmatrix}\notag \\  \notag\\
=\ &\begin{pmatrix}
[c_k]_q^{\flat}E^{\sharp}_{k-1}(c_1,\ldots,c_{k-1})_q-q^{c_{k-1}-1}E^{\sharp}_{k-2}(c_1,\ldots,c_{k-2})_q   \\
[c_k]_q^{\flat}E^{\sharp}_{k-2}(c_2,\ldots,c_{k-1})_q-q^{c_{k-1}-1}E^{\sharp}_{k-3}(c_2,\ldots,c_{k-2})_q \\
\end{pmatrix}\notag \\  \notag\\
=\ &\begin{pmatrix}
E^{\flat}_{k}(c_1,\ldots,c_{k})_q   \\
E^{\flat}_{k-1}(c_2,\ldots,c_{k})_q \\
\end{pmatrix}.\notag
\end{align}
Thus, by expanding the determinant \eqref{E-flat}, we can infer that 
\begin{align}
\displaystyle
\frac{E^{\flat}_{k}(c_1,\ldots,c_{k})_q}{E^{\flat}_{k-1}(c_2,\ldots,c_{k})_q}
&=[c_1]^{\sharp}_q-\cfrac{q^{c_1-1}}{\cfrac{E^{\flat}_{k-1}(c_2,\ldots,c_{k})_q}{E^{\flat}_{k-2}(c_3,\ldots,c_{k})_q}}=[c_1]^{\sharp}_q-\cfrac{q^{c_1-1}}{ [c_2]^{\sharp}_q- \cfrac{q^{c_2-1}}{\cfrac{E^{\flat}_{k-2}(c_3,\ldots,c_{k})_q}{E^{\flat}_{k-3}(c_4,\ldots,c_{k})_q}}}\notag \notag\\
&=\cdots=[[c_1,\ldots,c_k]]^{\flat}_q. \notag
\end{align}
\end{proof}

\textbf{Proof of Theorem \ref{Thm1}}:\\
By Proposition $4.9$ in \cite{MO1}, it follows that
\[
\displaystyle
q^{a_2+a_4+\cdots+a_{2m}}\sigma_{1,q}^{-a_1}\sigma_{2,q}^{a_2}\sigma_{1,q}^{-a_3}\sigma_{2,q}^{a_4}\cdots \sigma_{1,q}^{-a_{2m-1}} \sigma_{1,q}^{a_{2m}}=\sigma_{1,q}^{-c_1}S_q\sigma_{1,q}^{-c_2}S_q\cdots \sigma_{1,q}^{-c_k}S_q\sigma_{1,q}^{-1}
\]
and hence
\begin{align}
\displaystyle
[a_1,\ldots,a_{2m}]^{\flat}_q&=\sigma_{1,q}^{-a_1}\sigma_{2,q}^{a_2}\sigma_{1,q}^{-a_3}\sigma_{2,q}^{a_4}\cdots \sigma_{1,q}^{-a_{2m-1}} \sigma_{1,q}^{a_{2m}}\left(\frac{1}{1-q}\right)
\notag \\  \notag\\
&=\sigma_{1,q}^{-c_1}S_q\sigma_{1,q}^{-c_2}S_q\cdots \sigma_{1,q}^{-c_k}S_q\sigma_{1,q}^{-1}\left(\frac{1}{1-q}\right)\notag \\  \notag\\
&=\sigma_{1,q}^{-c_1}S_q\sigma_{1,q}^{-c_2}S_q\cdots \sigma_{1,q}^{-c_k}S_q\left(\frac{1}{1-q}\right)
\notag \\  \notag\\
&=[[c_1,\ldots,c_k]]^{\flat}_q. \notag
\end{align}
\qed

Through the above arguments, we have
\[\displaystyle
\left[\displaystyle\frac{r}{s}\right]_q^{\flat}=\sigma_{1,q}^{-c_1}S_q\sigma_{1,q}^{-c_2}S_q\cdots \sigma_{1,q}^{-c_k}S_q\sigma_{1,q}^{-1}\left(\frac{1}{1-q}\right).
\]



\subsection{Basic properties of the numerator and denominator of left $q$-rational numbers}
Morier-Genoud and Ovsienko give the basic properties of the numerator and denominator of right $q$-deformed rationals as follows \cite{MO1}:\\
\indent
For $i=1,2,\ldots,k$, we have
\[
\mathcal{R}^{\sharp}_{k}(q)=\mathcal{R}^{\sharp}(q), \ \ \ \
\mathcal{R}^{\sharp}_{i+1}(q)=[c_{i+1}]^{\sharp}_q\mathcal{R}^{\sharp}_{i}(q)-q^{c_i-1}\mathcal{R}^{\sharp}_{i-1}(q), 
\]
\[
\mathcal{S}^{\sharp}_{k}(q)=\mathcal{S}^{\sharp}(q), \ \ \ \
\mathcal{S}^{\sharp}_{i+1}(q)=[c_{i+1}]^{\sharp}_q\mathcal{S}^{\sharp}_{i}(q)-q^{c_i-1}\mathcal{S}^{\sharp}_{i-1}(q),   
\]
where the initial data
\[
\displaystyle
\mathcal{R}^{\sharp}_{0}(q)=1,\ \ \ \ \mathcal{R}^{\sharp}_{1}(q)=[c_1]_q,\ \ \ \ \mathcal{S}^{\sharp}_{0}(q)=0,\ \ \ \ \mathcal{S}^{\sharp}_{1}(q)=1,
\] 
then it follows that
\[
\displaystyle
\frac{\mathcal{R}^{\sharp}_{i}(q)}{\mathcal{S}^{\sharp}_{i}(q)}=[[c_1,\ldots,c_i]]_q.
\]\\

\indent
Similarly, we have the corresponding property for the left $q$-rationals as follows. 
For $i=1,2,\ldots,k$, we have
\[
\mathcal{R}^{\flat}_{k}(q)=\mathcal{R}^{\flat}(q), \ \ \ \
\mathcal{R}^{\flat}_{i+1}(q)=[c_{i+1}]^{\flat}_q\mathcal{R}^{\sharp}_{i}(q)-q^{c_i-1}\mathcal{R}^{\sharp}_{i-1}(q),
\]
\[
\mathcal{S}^{\flat}_{k}(q)=\mathcal{S}^{\flat}(q), \ \ \ \
\mathcal{S}^{\flat}_{i+1}(q)=[c_{i+1}]^{\flat}_q\mathcal{S}^{\sharp}_{i}(q)-q^{c_i-1}\mathcal{S}^{\sharp}_{i-1}(q),
\]
where the initial data
\[
\displaystyle
\mathcal{R}^{\flat}_{0}(q)=1,\ \ \ \ \mathcal{R}^{\flat}_{1}(q)=[c_1]^{\flat}_q,\ \ \ \ \mathcal{S}^{\flat}_{0}(q)=1-q,\ \ \ \ \mathcal{S}^{\flat}_{1}(q)=1,
\] 
then it follows that
\[
\displaystyle
\frac{\mathcal{R}^{\flat}_{i}(q)}{\mathcal{S}^{\flat}_{i}(q)}=[[c_1,\ldots,c_i]]^{\flat}_q.
\]\\

\medskip


\section{$q$-deformed Farey sum and $q$-deformed Farey triangles}\label{Sec3}

In this section, we give formulas corresponding to the $q$-deformed Farey sum of the left $q$-deformed rational numbers. We use this formula to obtain a $q$-deformed Farey tessellation and weighted triangulation on the left $q$-deformed rational numbers. For convenience,  from this section onwards, we always assume that the rational numbers are non-negative.  The case of negative rational numbers can be considered by symmetry.


\subsection{$q$-deformed Farey sum of left and right $q$-rational numbers}
\noindent
We consider two non-negative irreducible fractions $\displaystyle\frac{r}{s}$ and $\displaystyle\frac{r^{\prime}}{s^{\prime}}$ (we always asumme that $\displaystyle \frac{1}{0}$ is an irreducible fraction), then we say $\displaystyle\frac{r}{s}$ , $\displaystyle\frac{r^{\prime}}{s^{\prime}}$ are \textit {Farey neighbors} if $\left\vert sr^{\prime}-rs^{\prime}\right\vert=~1$. Different from the ordinary sum of fractions, we denote the \textit {Farey sum} of $\displaystyle\frac{r}{s}$ and $\displaystyle\frac{r^{\prime}}{s^{\prime}}$ by 
\begin{equation}\label{Farey_sum}
\displaystyle
\frac{r}{s}\#\frac{r^{\prime}}{s^{\prime}}:=\frac{r+r^{\prime}}{s+s^{\prime}}.
\end{equation}
The $q$-deformed Farey sum of right $q$-deformed rational numbers has been introduced in \cite{MO1}.
\begin{thm}[Morier-Genoud and Ovsienko \cite{MO1}]\label{THM-F-S-R}
For a positive rational number $\displaystyle \alpha=[[c_1,\ldots,c_k]]$ which is the Farey sum of \\ 
\begin{equation}\label{BETA}
\displaystyle
\beta=\left\{    
\begin{array}{
    cl}
    \ [[c_1,\ldots,c_{l}-1]] &  \text{for}  \  \   \  \ c_k=c_{k-1}=\cdots=c_{l+1}=2, \ c_{l}>2 ,\ 1\leq l\leq k   \\ & \\
    \ [[1]]  &  \text{for}  \ \ k=1 , \ c_k=2 \\ & \\
     \ [[1,1]]  &  \text{for}  \ \ c_k=c_{k-1}=\ldots=c_{2}=2, \ c_1=1 \\
\end{array} \right.
\end{equation}
and 
\begin{equation}\label{GAMMA}
\displaystyle
\gamma=\left\{    
\begin{array}{
    cl}
    \ [[c_1,\ldots,c_{k-1}]] &  \text{for}  \  \   \  \ k\geq2   \\ & \\
    \ [[]]  &  \text{for}  \ \ k=1,\\
\end{array} \right.
\end{equation}
if we assume that
\[
\displaystyle \left[\alpha \right]^{\sharp}_q=\frac{\mathcal{R}^{\sharp}_{\alpha}(q)}{\mathcal{S}^{\sharp}_{\alpha}(q)}, \ \ \ \
\displaystyle \left[\beta \right]^{\sharp}_q=\frac{\mathcal{R}^{\sharp}_{\beta}(q)}{\mathcal{S}^{\sharp}_{\beta}(q)},
\ \ \ \
\displaystyle \left[\gamma \right]^{\sharp}_q=\frac{\mathcal{R}^{\sharp}_{\gamma}(q)}{\mathcal{S}^{\sharp}_{\gamma}(q)},
\]
then

\begin{equation}\label{R-F-S}
\displaystyle
\frac{\mathcal{R}^{\sharp}_{\alpha}(q)}{\mathcal{S}^{\sharp}_{\alpha}(q)}
=\displaystyle\frac{\mathcal{R}^{\sharp}_{\beta}(q)+q^{c_k-1}\mathcal{R}^{\sharp}_{\gamma}(q)}{\mathcal{S}^{\sharp}_{\beta}(q)+q^{c_k-1}\mathcal{S}^{\sharp}_{\gamma}(q)}.
\end{equation}
\end{thm}

Hence, we define the $q$-defomed Farey sum $\#^{\sharp}_q$ of $\left[\beta \right]^{\sharp}_q$ and $\left[\gamma \right]^{\sharp}_q$ by \[\left[\beta \right]^{\sharp}_q \#^{\sharp}_q  \left[\gamma \right]^{\sharp}_q=\displaystyle\frac{\mathcal{R}^{\sharp}_{\beta}(q)+q^{c_k-1}\mathcal{R}^{\sharp}_{\gamma}(q)}{\mathcal{S}^{\sharp}_{\beta}(q)+q^{c_k-1}\mathcal{S}^{\sharp}_{\gamma}(q)}.\]

\medskip

\begin{exam}
\ \\
{\normalfont
\indent
$\displaystyle\frac{12}{5}=[[3,2,3]]$, $\displaystyle\frac{7}{3}=[[3,2,2]]$, $\displaystyle\frac{5}{2}=[[3,2]]$,  then
\[
\displaystyle
\left[\frac{12}{5}\right]^{\sharp}_q=\frac {1+2\,q+3\,{q}^{2}+3\,{q}^{3}+2\,{q}^{4}+{q}^{5}}{1+q+2\,{q}^{2
}+{q}^{3}}
,
\]
\[
\displaystyle
\left[\frac{7}{3}\right]^{\sharp}_q=\frac {1+2\,q+2\,{q}^{2}+{q}^{3}+{q}^{4}}{1+q+{q}^{2}},
\ \ \ \ \
\left[\frac{5}{2}\right]^{\sharp}_q=\frac {1+2\,q+{q}^{2}+{q}^{3}}{1+q},
\]
\[
\displaystyle
\left[\frac{7}{3}\right]^{\sharp}_q \#^{\sharp}_q \left[\frac{5}{2}\right]^{\sharp}_q= \frac{(1+2\,q+2\,{q}^{2}+{q}^{3}+{q}^{4})+q^2(1+2\,q+{q}^{2}+{q}^{3})}{(1+q+{q}^{2})+q^2(1+q)}=\left[\frac{12}{5}\right]^{\sharp}_q.
\]
}
\end{exam}

\medskip
Similarly, we consider a left $q$-deformed rational number for a left $q$-deformed Farey sum. The following theorem gives the formula for the $q$-deformed Farey sum of a left $q$-deformed rational number. It is interesting to note that formula \eqref{L-F-S} forms a formal symmetry with the formula \eqref{R-F-S}.

\begin{thm}\label{THM-F-S-L}
For a rational number $\displaystyle \alpha=[[c_1,\ldots,c_k]]$ which is the Farey sum of $\beta$ and $\gamma$ defined by \eqref{BETA} and \eqref{GAMMA}, if we assume that
\[
\displaystyle \left[\alpha \right]^{\flat}_q=\frac{\mathcal{R}_{\alpha}^{\flat}(q)}{\mathcal{S}_{\alpha}^{\flat}(q)}, \ \ \ \
\displaystyle \left[\beta \right]^{\flat}_q=\frac{\mathcal{R}_{\beta}^{\flat}(q)}{\mathcal{S}_{\beta}^{\flat}(q)},
\ \ \ \
\displaystyle \left[\gamma \right]^{\flat}_q=\frac{\mathcal{R}_{\gamma}^{\flat}(q)}{\mathcal{S}_{\gamma}^{\flat}(q)},
\]
 then

\begin{equation}\label{L-F-S}
\displaystyle
\frac{\mathcal{R}_{\alpha}^{\flat}(q)}{\mathcal{S}_{\alpha}^{\flat}(q)} 
=\displaystyle\frac{q^{k-l+1}\mathcal{R}_{\beta}^{\flat}(q)+\mathcal{R}_{\gamma}^{\flat}(q)}{q^{k-l+1}\mathcal{S}_{\beta}^{\flat}(q)+\mathcal{S}_{\gamma}^{\flat}(q)},
\end{equation}
where $c_k=c_{k-1}=\cdots=c_{l+1}=2, \ c_{l}>2 ,\ 1\leq l\leq k$.\\
In particular, for $k=1$, we have 
\[
\displaystyle
\frac{\mathcal{R}_{\alpha}^{\flat}(q)}{\mathcal{S}_{\alpha}^{\flat}(q)} 
=\displaystyle\frac{q\mathcal{R}_{\beta}^{\flat}(q)+\mathcal{R}_{\gamma}^{\flat}(q)}{q\mathcal{S}_{\beta}^{\flat}(q)+\mathcal{S}_{\gamma}^{\flat}(q)}.
\]
\end{thm}
Hence, we define the $q$-defomed Farey sum $\#^{\flat}_q$ of $\left[\beta \right]^{\flat}_q$ and $\left[\gamma \right]^{\flat}_q$ by \[\left[\beta \right]^{\flat}_q \#^{\flat}_q  \left[\gamma \right]^{\flat}_q=\displaystyle\frac{q^{k-l+1}\mathcal{R}_{\beta}^{\flat}(q)+\mathcal{R}_{\gamma}^{\flat}(q)}{q^{k-l+1}\mathcal{S}_{\beta}^{\flat}(q)+\mathcal{S}_{\gamma}^{\flat}(q)}.\]

\medskip

\begin{proof}
Suppose that $\alpha=[[c_1,\dots,c_l,2^{(k-l)}]]$, where $2^{(k-l)}$ stands for $k-l$ copies of $2$, $\beta=[[c_1,\dots,c_l-1]]$, $\gamma=[[c_1,\dots,c_l,2^{(k-l-1)}]]$, then

\begin{align}
\displaystyle
\mathcal{R}_{\alpha}^{\flat}(q)&=E^{\flat}_{k}(c_1,\dots,c_l,2^{(k-l)})_q \notag \\
&=[2]^{\flat}_qE^{\sharp}_{k-1}(c_1,\dots,c_l,2^{(k-l-1)})_q-qE^{\sharp}_{k-2}(c_1,\dots,c_l,2^{(k-l-2)})_q \notag \\
&=(1+q^2+q^3)E^{\sharp}_{k-2}(c_1,\dots,c_l,2^{(k-l-2)})_q-(q+q^3)E^{\sharp}_{k-3}(c_1,\dots,c_l,2^{(k-l-3)})_q, \notag
\end{align}
and
\begin{align}
\displaystyle
\mathcal{R}_{\gamma}^{\flat}(q)&=E^{\flat}_{k-1}(c_1,\dots,c_l,2^{(k-l-1)})_q \notag \\
&=[2]^{\flat}_qE^{\sharp}_{k-2}(c_1,\dots,c_l,2^{(k-l-2)})_q-qE^{\sharp}_{k-3}(c_1,\dots,c_l,2^{(k-l-3)})_q\notag \\
&=(1+q^2)E^{\sharp}_{k-2}(c_1,\dots,c_l,2^{(k-l-2)})_q-qE^{\sharp}_{k-3}(c_1,\dots,c_l,2^{(k-l-3)})_q.\notag 
\end{align}

Thus, by \eqref{E-flat-2},
\begin{align}
\displaystyle
\mathcal{R}^{\flat}_{\alpha}(q)-\mathcal{R}^{\flat}_{\gamma}(q)&=E^{\sharp}_{k+1}(c_1,\dots,c_l,2^{(k-l+1)})_q-E^{\sharp}_{k}(c_1,\dots,c_l,2^{(k-l)})_q. \notag
\end{align}

On the other hand,
\begin{align}
\displaystyle
\mathcal{R}^{\flat}_{\beta}(q)&=E^{\flat}_{l}(c_1,\dots,c_l-1)_q\notag \\
&=[c_l-1]^{\flat}_qE^{\sharp}_{l-1}(c_1,\dots,c_{l-1})_q-q^{c_{l-1}-1}E^{\sharp}_{l-2}(c_1,\dots,c_{l-2})_q\notag \\
&=([c_l]_q-q^{c_l-2})E^{\sharp}_{l-1}(c_1,\dots,c_{l-1})_q-q^{c_l-2}E^{\sharp}_{l-2}(c_1,\dots,c_{l-2})_q\notag \\
&=E^{\sharp}_{l}(c_1,\dots,c_{l})_q-q^{c_l-2}E^{\sharp}_{l-1}(c_1,\dots,c_{l-1})_q.\notag
\end{align}

Again, by \eqref{E-flat-2}, we have
\begin{align}
\displaystyle
q^{k-l+1}\mathcal{R}^{\flat}_{\beta}(q)&=q^{k-l}(q(E^{\sharp}_{l}(c_1,\dots,c_{l})_q-q^{c_l-2}E^{\sharp}_{l-1}(c_1,\dots,c_{l-1})_q))\notag \\
&=q^{k-l}(E^{\sharp}_{l+1}(c_1,\dots,c_{l},2)_q-E^{\sharp}_{l}(c_1,\dots,c_{l})_q)\notag \\
&=E^{\sharp}_{k+1}(c_1,\dots,c_l,2^{(k-l+1)})_q-E^{\sharp}_{k}(c_1,\dots,c_l,2^{(k-l)})_q. \notag
\end{align}

Hence, we proved that,
\[
\displaystyle
\mathcal{R}_{\alpha}^{\flat}(q)=q^{k-l+1}\mathcal{R}_{\beta}^{\flat}(q)+\mathcal{R}_{\gamma}^{\flat}(q).
\]
The proof of $\mathcal{S}_{\alpha}^{\flat}(q)=q^{k-l+1}\mathcal{S}_{\beta}^{\flat}(q)+\mathcal{S}_{\gamma}^{\flat}(q)$ is similar.

\end{proof}

\medskip

\begin{exam}
\ \\
{\normalfont
\indent
(1) Since $\displaystyle\frac{12}{5}=[[3,2,3]]$, $\displaystyle\frac{7}{3}=[[3,2,2]]$, $\displaystyle\frac{5}{2}=[[3,2]]$, then by Theorem \ref{THM-F-S-L}, it follows that
\[
\displaystyle
\left[\frac{12}{5}\right]^{\flat}_q=\frac {1+2\,q+2\,{q}^{2}+2\,{q}^{3}+3\,{q}^{4}+{q}^{5}+{q}^{6}}{1+q+{
q}^{2}+{q}^{3}+{q}^{4}},
\]
\[
\displaystyle
\left[\frac{7}{3}\right]^{\flat}_q=\frac{1+q+q^2+2q^3+q^4+q^5}{1+q^2+q^3},
\ \ \ \ \
\left[\frac{5}{2}\right]^{\flat}_q=\frac{1+q+q^2+q^3+q^4}{1+q^2},
\]
\[
\displaystyle
\left[\frac{7}{3}\right]^{\flat}_q \#^{\flat}_q  \left[\frac{5}{2}\right]^{\flat}_q= \frac{q(1+q+q^2+2q^3+q^4+q^5)+(1+q+q^2+q^3+q^4)}{q(1+q^2+q^3)+(1+q^2)}=\left[\frac{12}{5}\right]^{\flat}_q.
\]
\\

(2) Since $\displaystyle\frac{7}{2}=[[4,2]]$, $\displaystyle\frac{3}{1}=[[3]]$, $\displaystyle\frac{4}{1}=[[4]]$, then by Theorem \ref{THM-F-S-L}, it follows that
\[
\displaystyle
\left[\frac{7}{2}\right]^{\flat}_q=\frac {1+q+2\,{q}^{2}+{q}^{3}+{q}^{4}+{q}^{5}}{1+{q}^{2}},
\\
\displaystyle
\left[\frac{3}{1}\right]^{\flat}_q=\frac{1+q+q^3}{1},
\\
\left[\frac{4}{1}\right]^{\flat}_q=\frac{1+q+q^2+q^4}{1},
\\
\]
\[
\displaystyle
\left[\frac{3}{1}\right]^{\flat}_q \#^{\flat}_q  \left[\frac{4}{1}\right]^{\flat}_q=\frac{q^2(1+q+q^3)+(1+q+q^2+q^4)}{q^2+1} =\left[\frac{7}{2}\right]^{\flat}_q.
\]

(3) Since $\displaystyle\frac{9}{4}=[[3,2,2,2]]$, $\displaystyle\frac{2}{1}=[[2]]$, $\displaystyle\frac{7}{3}=[[3,2,2]]$, then by Theorem \ref{THM-F-S-L}, it follows that
\[
\displaystyle
\left[\frac{9}{4}\right]^{\flat}_q=\frac {1+q+{q}^{2}+2\,{q}^{3}+2\,{q}^{4}+{q}^{5}+{q}^{6}}{1+{q}^{2}+{
q}^{3}+{q}^{4}},
\]
\[\displaystyle
\left[\frac{7}{3}\right]^{\flat}_q=\frac {1+q+{q}^{2}+2\,{q}^{3}+{q}^{4}+{q}^{5}}{1+{q}^{2}+{q}^{3}},
\\
\left[\frac{2}{1}\right]^{\flat}_q=\frac{1+q^2}{1},
\]
\[\displaystyle
\left[\frac{2}{1}\right]^{\flat}_q \#^{\flat}_q  \left[\frac{7}{3}\right]^{\flat}_q= \frac{q^4(1+q^2)+(1+q+{q}^{2}+2\,{q}^{3}+{q}^{4}+{q}^{5})}{q^4+(1+{q}^{2}+{q}^{3})}=\left[\frac{9}{4}\right]^{\flat}_q.
\]

}
\end{exam}

\subsection{$q$-deformed Farey tessellation about left $q$-deformed rational numbers}

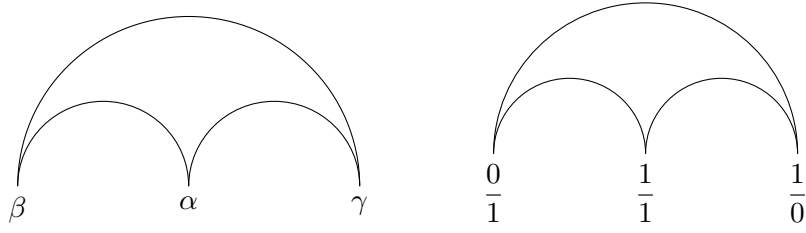
\begin{figure}
\begin{center}
\begin{tikzpicture}[scale=4.5]
    
   \draw (1,0) arc (0:180:.5);

    \draw  (0,0)  node[below]{$\displaystyle\beta$};
    \draw  (1,0)  node[below]{$\displaystyle\gamma$};

    \draw (1,0) arc (0:180:.25);
    \draw  (.5,0) node[below]{$\displaystyle\alpha$};
    \draw (.5,0) arc (0:180:.25);
    

\end{tikzpicture} \hspace{1cm}
\begin{tikzpicture}[scale=4]
    
    \draw (1,0) arc (0:180:.5);

    \draw  (0,0)  node[below]{$\displaystyle\frac{0}{1}$};
    \draw  (1,0)  node[below]{$\displaystyle\frac{1}{0}$};

    \draw (1,0) arc (0:180:.25);
    \draw  (.5,0) node[below]{$\displaystyle\frac{1}{1}$};
    \draw (.5,0) arc (0:180:.25);
    
\end{tikzpicture}

\end{center}
\caption{Farey triangle (left) and the initial Farey triangle (right).}
    \label{Fig.1}
\end{figure}

In this section, following \cite{MO1, BBL}, we discuss a relationship between left $q$-deformed rational numbers and Farey tessellation (see \cite{HW} for more details). We assume that all rational numbers are represented as irreducible fractions. We order the elements of $\mathbb{Q}_{>0}\cup\{\infty\}$ by horizontal segment drawn in the plane, then Farey tessellation consists of all triangles whoose forms are as in the left of Figure \ref{Fig.1} (a rational number $\displaystyle \alpha$ which is the Farey sum of $\beta$ and $\gamma$ defined by \eqref{BETA} and \eqref{GAMMA}), and each vertex corresponds to a rational number, and any two vertices that are Farey neighbors are connected by a semicircle. We call these triangles Farey triangles, and the initial Farey triangle is on the right of Figure \ref{Fig.1}.

\medskip

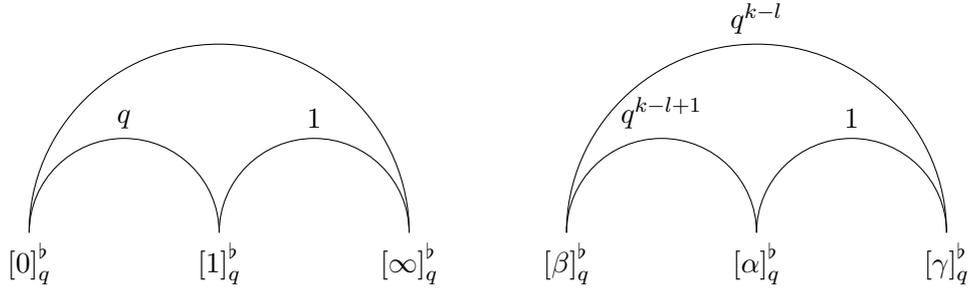
\begin{figure}
\begin{center}
\begin{tikzpicture}[scale=5]
    
    \draw (1,0) arc (0:180:.5);

    \draw  (0,0)  node[below]{$\displaystyle\left[0\right]_q^{\flat}$};
    \draw  (1,0)  node[below]{$\displaystyle\left[\infty\right]_q^{\flat}$};

    \draw (1,0) arc (0:180:.25);
    \draw  (.5,0) node[below]{$\displaystyle\left[1\right]_q^{\flat}$};
    \draw (.5,0) arc (0:180:.25);
    
\draw (.25,.25) node[above]{$q$};
\draw (.75,.25) node[above]{$1$};
\end{tikzpicture}\hspace{1cm}
\begin{tikzpicture}[scale=5]
    
   \draw (1,0) arc (0:180:.5);

    \draw  (0,0)  node[below]{$\displaystyle\left[\beta\right]_q^{\flat}$};
    \draw  (1,0)  node[below]{$\displaystyle\left[\gamma\right]_q^{\flat}$};

    \draw (1,0) arc (0:180:.25);
    \draw  (.5,0) node[below]{$\displaystyle\left[\alpha\right]_q^{\flat}$};
    \draw (.5,0) arc (0:180:.25);
    
\draw (.25,.25) node[above]{$q^{k-l+1}$};
\draw (.75,.25) node[above]{$1$};
\draw (.5,.5) node[above]{$q^{k-l}$};

\end{tikzpicture}

\end{center}
\caption{The initial $q$-deformed Farey triangle (left) and the $q$-deformed Farey triangle (right).}
    \label{Fig.2}
\end{figure}

\medskip

\begin{figure}

\begin{center}
\begin{tikzpicture}[scale=14]
   
    \draw (1,0) arc (0:180:.5);

    \draw (0,0) node[below]{$\displaystyle\left[\frac{0}{1}\right]^{\flat}_q$};
    \draw (1,0) node[below]{$\displaystyle\left[\frac{1}{0}\right]^{\flat}_q$};

    \draw (1,0) arc (0:180:.25);
    \draw (.5,0) node[below]{$\displaystyle\left[\frac{1}{1}\right]^{\flat}_q$};
    \draw (.5,0) arc (0:180:.25);

    \draw (1,0) arc (0:180:1/6);
    \draw  (2/3,0)  node[below]{$\displaystyle\left[\frac{2}{1}\right]^{\flat}_q$};
    \draw (1/3,0) arc (0:180:1/12);
    \draw (1/6,0) arc (0:180:1/12);
    \draw (1/6,0) node[below]{$\displaystyle\left[\frac{1}{3}\right]^{\flat}_q$};
    
    \draw (2/3,0) arc (0:180:1/12);

    \draw (1/3,0) arc (0:180:1/6);
    \draw (1/3,0) node[below]{$\displaystyle\left[\frac{1}{2}\right]^{\flat}_q$};
    \draw (1/2,0) arc (0:180:1/12);
    \draw (5/6,0) node[below]{$\displaystyle\left[\frac{3}{1}\right]^{\flat}_q$};
    \draw (5/6,0) arc (0:180:1/12);
    \draw (1,0) arc (0:180:1/12);
    
    \draw (1/4,1/4)node[above]{$q$};
    \draw (3/4,1/4)node[above]{$1$};
    
\draw (1/12,1/12) node[above]{$q^3$};
\draw (3/12,1/12) node[above]{$1$};

\draw (1/6,1/6) node[above]{$q^2$};
\draw (5/12,1/12) node[above]{$1$};

\draw (7/12,1/12) node[above]{$q$};
\draw (10/12,1/6) node[above]{$1$};

\draw (11/12,1/12) node[above]{$1$};
\draw (9/12,1/12) node[above]{$q$};

    \end{tikzpicture}
    \end{center}
    \caption{A part of the Farey tessellation with weights carried by the edges and left $q$-deformed rational numbers labeling the vertices.}
    \label{Fig.3}
\end{figure}
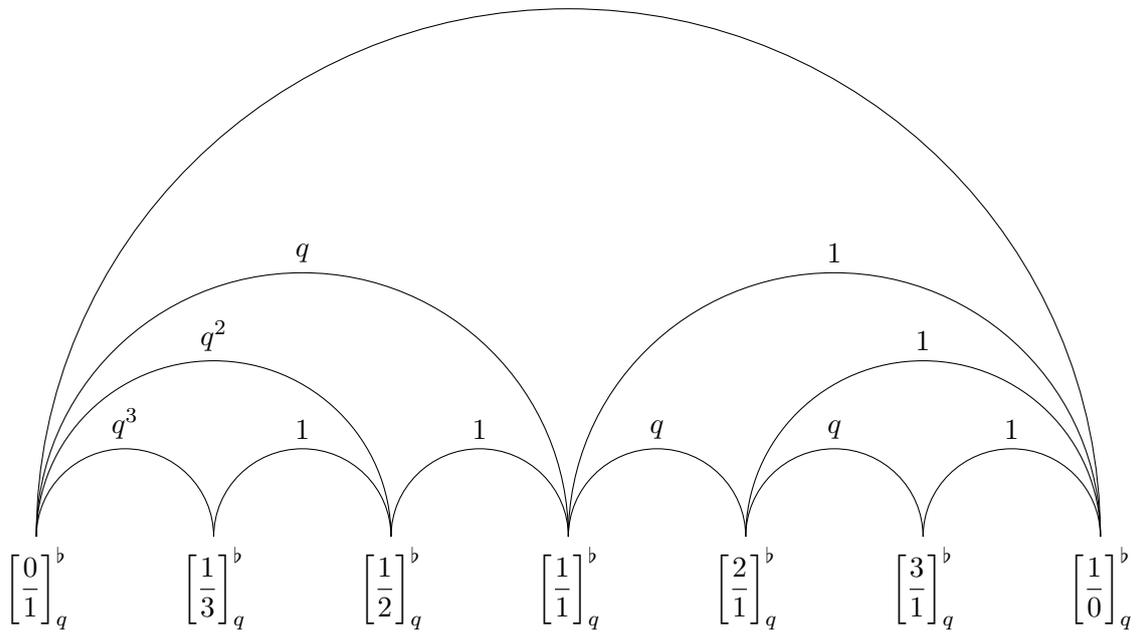

The $q$-deformed Farey tessellation considered in \cite{MO1} and \cite{BBL} is composed of $q$-deformed Farey triangles which are obtained by basing on the original Farey triangle and each vertex is a right $q$-deformed rational number and each edge is weighted. Every $q$-deformed Farey triangle can be obtained according to the laws of Theorem \ref{THM-F-S-R}. Now we replace Theorem \ref{THM-F-S-R} with Theorem \ref{THM-F-S-L}, and setting the initial $q$-deformed Farey triangle as the left of Figure \ref{Fig.2}, then we can obtain a new Farey tessellation consisting of $q$-deformed Farey triangles as in the right of Figure \ref{Fig.2}. Each vertex of a $q$-deformed Farey triangle corresponds to a left $q$-deformed rational number (as a simple example, see Figure \ref{Fig.3}).

\medskip

Following \cite[Section 2.2]{BBL}, we choose two infinitely close Farey triangle sequences from the left and right sides near the rational number $\displaystyle\alpha$. Considering the Farey tessellation according to the laws of Theorem \ref{THM-F-S-L}, then we find that the Farey triangle sequence on the left side of $\displaystyle\alpha$ converges to exactly one point. However, when $q$ is not equal to $1$, the one on the right side of $\displaystyle\alpha$ cannot converge to a point. Thus we obtain a figure with mirror symmetry to \cite[Figure~5]{BBL}.

\subsection{Weighted triangulation about left $q$-deformed rational numbers}

\noindent
Consider a positive rational number $\displaystyle \alpha=[a_1,\ldots,a_{2m}]=[[c_1,\ldots,c_k]]$ which is the Farey sum of $\beta$ and $\gamma$ defined by \eqref{BETA} and \eqref{GAMMA}. According to \cite{MO3}, it follows that $\displaystyle\alpha$ corresponds to a triangulation as in Figure \ref{Fig.4}. If we give the initial values as in Figure \ref{Fig.5}, then the remaining vertices can be computed according to the Farey sum.  

\begin{figure}[H]
\begin{center}
\begin{tikzpicture}[scale=4]

\draw  (0,0)  -- (.5,.6) --(1,0)--(0,0);
\draw   (.5,.6) --(.25,0);
\draw   (.5,.6) --(.75,0);
\draw (.5,.24) node[above]{$\cdots$};

\draw  (1,0)  -- (1.5,.6) --(.5,.6);
\draw   (1,0) --(0.75,.6);
\draw   (1,0) --(1.25,.6);
\draw (1,.24) node[above]{$\cdots$};

\draw(1,0)--(1.5,0);

\draw[dashed](1.5,0)--(2,0);
\draw[dashed](1.5,.3)--(2,.3);
\draw[dashed](1.5,.6)--(2,.6);

\draw  (2,.6)--(2.5,.6)--(2,0);

\draw  (2,0)--(3,0)--(2.5,.6);
\draw   (2.5,.6) --(2.25,0);
\draw   (2.5,.6) --(2.75,0);
\draw (2.5,.24) node[above]{$\cdots$};

\draw  (2.5,.6)--(3.5,.6)--(3,0);
\draw   (3,0) --(2.75,.6);
\draw   (3,0) --(3.25,.6);
\draw (3,.24) node[above]{$\cdots$};

\draw (.5,-.1) node{$a_1$};
\draw [->] (.45,-.1)--(0,-.1);
\draw [->] (.55,-.1)--(1,-.1);

\draw (1.5,-.1) node{$a_3$};
\draw [->] (1.45,-.1)--(1,-.1);

\draw (2.5,-.1) node{$a_{2m-1}$};
\draw [->] (2.3,-.1)--(2,-.1);
\draw [->] (2.70,-.1)--(3,-.1);

\draw (1,0.7) node{$a_2$};
\draw [->] (.95,0.7)--(.5,0.7);
\draw [->] (1.05,0.7)--(1.5,0.7);

\draw (2,0.7) node{$a_{2m-2}$};
\draw [->] (2.15,0.7)--(2.5,0.7);

\draw (3,0.7) node{$a_{2m}$};
\draw [->] (2.9,0.7)--(2.5,0.7);
\draw [->] (3.1,0.7)--(3.5,0.7);


    \end{tikzpicture}
    \end{center}
    \caption{Triangulation of $\displaystyle\alpha$.}
    
    \label{Fig.4}
\end{figure}

\begin{figure}[H]
\begin{center}
\begin{tikzpicture}[scale=4]

    \draw  (0,0) node[below]{$\displaystyle\frac{0}{1}$} -- (.5,.6) node[above]{$\displaystyle\frac{1}{0}$}--(.5,0)node[below]{$\displaystyle\frac{1}{1}$}--(0,0);

    \draw  (1,0) node[below]{$\displaystyle\frac{0}{1}$} -- (1,0.6) node[above]{$\displaystyle\frac{1}{0}$}--(1.5,.6)node[above]{$\displaystyle\frac{1}{1}$}--(1,0);

    \end{tikzpicture}
    \end{center}
    \caption{Initial settings of triangulations for the cases of $\displaystyle\alpha>1$(left) and $0<\displaystyle\alpha\leq1$(right).}
    
    \label{Fig.5}
\end{figure}
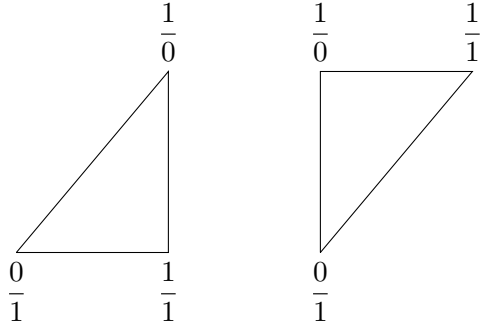

\begin{exam}
{\normalfont

The triangulations of $\displaystyle \frac{8}{11}=[0,1,2,1,1,1]$ and $\displaystyle \frac{11}{8}=[1,2,2,1]$ can be expressed as Figure \ref{Fig.6}. 

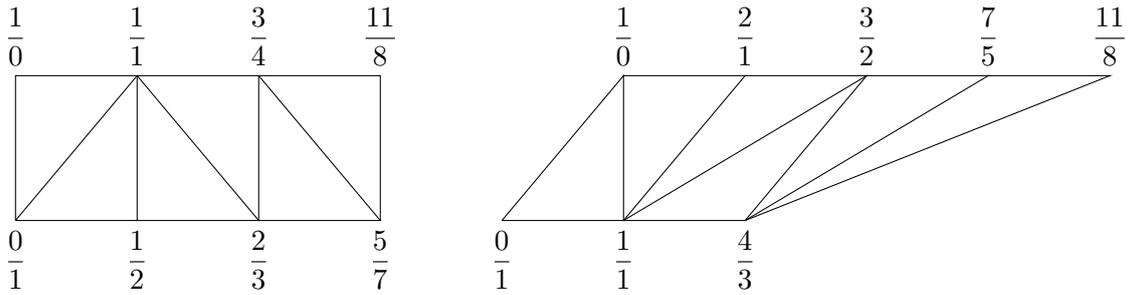
\begin{figure}[H]
\begin{center}
\begin{tikzpicture}[scale=3.2]

    \draw  (0,0) node[below]{$\displaystyle\frac{0}{1}$} -- (0,0.6) node[above]{$\displaystyle\frac{1}{0}$}--(.5,.6)node[above]{$\displaystyle\frac{1}{1}$}--(0,0);

\draw  (0,0)--(.5,0)node[below]{$\displaystyle\frac{1}{2}$}--(.5,.6);
\draw  (.5,0)--(1,0)node[below]{$\displaystyle\frac{2}{3}$}--(.5,.6);
\draw  (1,0)--(1,.6)node[above]{$\displaystyle\frac{3}{4}$}--(.5,.6);
\draw  (1,.6)--(1.5,0)node[below]{$\displaystyle\frac{5}{7}$}--(1,0);
\draw  (1.5,0)--(1.5,.6)node[above]{$\displaystyle\frac{11}{8}$}--(1,.6);


\draw  (2,0) node[below]{$\displaystyle\frac{0}{1}$} -- (2.5,.6) node[above]{$\displaystyle\frac{1}{0}$}--(2.5,0)node[below]{$\displaystyle\frac{1}{1}$}--(2,0);

\draw  (2.5,0)--(3,.6)node[above]{$\displaystyle\frac{2}{1}$}--(2.5,.6);

\draw  (2.5,0)--(3.5,.6)node[above]{$\displaystyle\frac{3}{2}$}--(3,.6);
\draw  (2.5,0)--(3,0)node[below]{$\displaystyle\frac{4}{3}$}--(3.5,.6);
\draw  (3,0)--(4,.6)node[above]{$\displaystyle\frac{7}{5}$}--(3.5,.6);
\draw  (3,0)--(4.5,.6)node[above]{$\displaystyle\frac{11}{8}$}--(4,.6);


    \end{tikzpicture}
    \end{center}
    \caption{triangulations of $\displaystyle\frac{8}{11}$ (left) and $\displaystyle\frac{11}{8}$ (right).}
    
    \label{Fig.6}
\end{figure}

}
\end{exam}

\medskip

\begin{figure}[H]
\begin{center}
\begin{tikzpicture}[scale=4.5]

    \draw  (0,0) node[below]{$\displaystyle\left[\beta\right]_q^{\sharp}$} -- (.5,.6) node[above]{$\displaystyle\left[\alpha\right]_q^{\sharp}$}--(.7,0)node[below]{$\displaystyle\left[\gamma\right]_q^{\sharp}$}--(0,0);

\draw(.35,0)node[below]{$\displaystyle q^{c_k-2}$};
\draw(.2,.3)node[above]{$\displaystyle 1$};
\draw(.6,.3)node[right]{$\displaystyle q^{c_k-1}$};

    \end{tikzpicture}
    \end{center}
    \caption{The triangle set by Theorem \ref{THM-F-S-R}.}
    \label{Fig.7-1}
\end{figure}
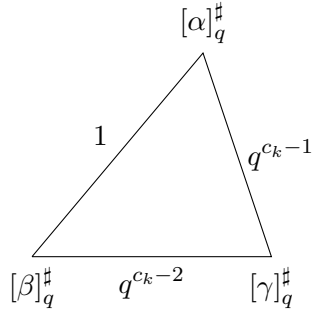

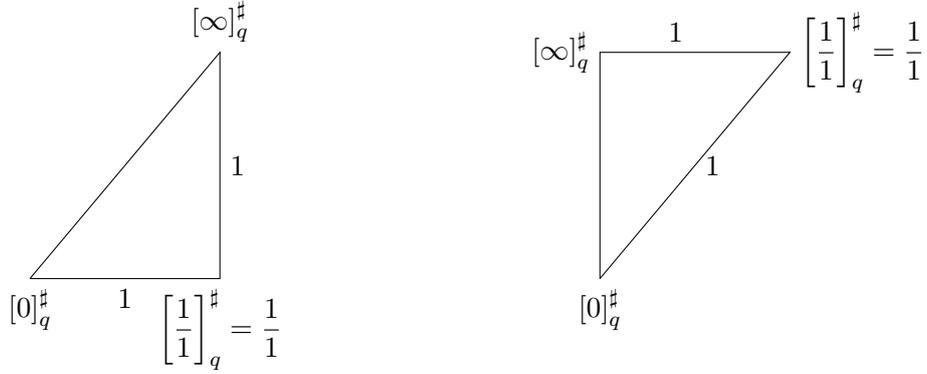
\begin{figure}[H]
\begin{center}
\begin{tikzpicture}[scale=5]

    \draw  (0,0) node[below]{$\displaystyle\left[0\right]_q^{\sharp}$} -- (.5,.6) node[above]{$\displaystyle\left[\infty\right]_q^{\sharp}$}--(.5,0)node[below]{$\displaystyle\left[\frac{1}{1}\right]_q^{\sharp}=\frac{1}{1}$}--(0,0);

    \draw  (1.5,0) node[below]{$\displaystyle\left[0\right]_q^{\sharp}$} -- (1.5,0.6) node[left]{$\displaystyle\left[\infty\right]_q^{\sharp}$}--(2,.6)node[right]{$\displaystyle\left[\frac{1}{1}\right]_q^{\sharp}=\frac{1}{1}$}--(1.5,0);

\draw(.25,0)node[below]{$\displaystyle 1$};
\draw(.5,.3)node[right]{$\displaystyle1$};

\draw(1.7,.6)node[above]{$\displaystyle 1$};
\draw(1.75,.3)node[right]{$\displaystyle 1$};

    \end{tikzpicture}
    \end{center}
    \caption{The initial settings of our's weighted triangulation for the cases $\displaystyle\alpha>1$(left) and $0<\displaystyle\alpha\leq1$(right).}
    \label{Fig.8-1}
\end{figure}

\indent
Now, let us consider the $q$-deformation of the triangulation which is called weighted triangulation.  For the vertices and edges of the two kinds of triangles in the triangulation, we will set them with the $q$-deformed Farey sum from Theorem \ref{THM-F-S-R} in Figure \ref{Fig.7-1}, and the initial setting is as Figure \ref{Fig.8-1} (See \cite{MO1} for details).

\indent
Similary, for the vertices and edges of the two kinds of triangles in the triangulation, we will set them with the $q$-deformed Farey sum from Theorem \ref{THM-F-S-L} in Figure \ref{Fig.7}, and the initial setting is as Figure \ref{Fig.8}. Different from the above case, the weights of the edges with  endpoints $[\alpha]_q$ and $[\gamma]_q$ become $1$, and the weights of the edges with endpoints $[\alpha]_q$ and $[\beta]_q$ become some power of $q$, which is symmetric to the case of the right $q$-deformed rational numbers. This phenomenon also corresponding to the fact that for the $q$-deformed Farey sum, there is a symmetric between the case of the right $q$-deformed  rational numbers and the case of the left $q$-deformed rational numbers.

\begin{figure}[H]
\begin{center}
\begin{tikzpicture}[scale=4.5]

    \draw  (0,0) node[below]{$\displaystyle\left[\beta\right]_q^{\flat}$} -- (.5,.6) node[above]{$\displaystyle\left[\alpha\right]_q^{\flat}$}--(.7,0)node[below]{$\displaystyle\left[\gamma\right]_q^{\flat}$}--(0,0);

\draw(.35,0)node[below]{$\displaystyle q^{k-l}$};
\draw(.2,.3)node[above]{$\displaystyle q^{k-l+1}$};
\draw(.6,.3)node[right]{$\displaystyle1$};

    \end{tikzpicture}
    \end{center}
    \caption{The triangle set by Theorem \ref{THM-F-S-L}.}
    \label{Fig.7}
\end{figure}
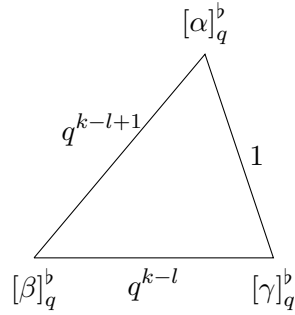

\begin{figure}[H]
\begin{center}
\begin{tikzpicture}[scale=5]

    \draw  (0,0) node[below]{$\displaystyle\left[0\right]_q^{\flat}$} -- (.5,.6) node[above]{$\displaystyle\left[\infty\right]_q^{\flat}$}--(.5,0)node[below]{$\displaystyle\left[\frac{1}{1}\right]_q^{\flat}=\frac{q}{1}$}--(0,0);

    \draw  (1.5,0) node[below]{$\displaystyle\left[0\right]_q^{\flat}$} -- (1.5,0.6) node[left]{$\displaystyle\left[\infty\right]_q^{\flat}$}--(2,.6)node[right]{$\displaystyle\left[\frac{1}{1}\right]_q^{\flat}=\frac{q}{1}$}--(1.5,0);

\draw(.25,0)node[below]{$\displaystyle q$};
\draw(.5,.3)node[right]{$\displaystyle1$};

\draw(1.7,.6)node[above]{$\displaystyle 1$};
\draw(1.75,.3)node[right]{$\displaystyle q$};

    \end{tikzpicture}
    \end{center}
    \caption{The initial settings of our's weighted triangulation for the cases $\displaystyle\alpha>1$(left) and $0<\displaystyle\alpha\leq1$(right).}
    \label{Fig.8}
\end{figure}

\medskip

\medskip


\begin{exam}
{\normalfont
The weighted triangulations of $\displaystyle\left[\frac{8}{11}\right]_q^{\sharp}=[0,1,2,1,1,1]_q^{\sharp}=[[1,4,3]]_q^{\sharp}$ and $\displaystyle\left[\frac{11}{8}\right]_q^{\sharp}=[1,2,2,1]^{\sharp}_q=[[2,2,4]]_q^{\sharp}$ can be expressed as Figures \ref{Fig.9-1} and \ref{Fig.10-1}. 
\begin{figure}[H]
\begin{center}
\begin{tikzpicture}[scale=5.5]

    \draw  (0,0) node[below]{$\displaystyle\left[0\right]_q^{\sharp}$} -- (0,0.6) node[above]{$\displaystyle\left[\infty\right]_q^{\sharp}$}--(.5,.6)node[above]{$\displaystyle\left[\frac{1}{1}\right]_q^{\sharp}$}--(0,0);

\draw  (0,0)--(.5,0)node[below]{$\displaystyle\left[\frac{1}{2}\right]_q^{\sharp}$}--(.5,.6);
\draw  (.5,0)--(1,0)node[below]{$\displaystyle\left[\frac{2}{3}\right]_q^{\sharp}$}--(.5,.6);
\draw  (1,0)--(1,.6)node[above]{$\displaystyle\left[\frac{3}{4}\right]_q^{\sharp}$}--(.5,.6);
\draw  (1,.6)--(1.5,0)node[below]{$\displaystyle\left[\frac{5}{7}\right]_q^{\sharp}$}--(1,0);
\draw  (1.5,0)--(1.5,.6)node[above]{$\displaystyle\left[\frac{11}{8}\right]_q^{\sharp}$}--(1,.6);

\draw(.25,0)node[above]{$\displaystyle 1$};
\draw(.25,.3)node[above]{$\displaystyle 1$};
\draw(.25,.6)node[above]{$\displaystyle 1$};
\draw(.5,.3)node[right]{$\displaystyle q$};
\draw(.75,.0)node[above]{$\displaystyle 1$};
\draw(.75,.3)node[above]{$\displaystyle q^2$};
\draw(.75,.6)node[above]{$\displaystyle q^3$};
\draw(1,.3)node[right]{$\displaystyle 1$};
\draw(1.25,.0)node[above]{$\displaystyle 1$};
\draw(1.25,.3)node[above]{$\displaystyle q$};
\draw(1.25,.6)node[above]{$\displaystyle q^2$};
\draw(1.5,.3)node[right]{$\displaystyle 1$};

    \end{tikzpicture}
    \end{center}
    \caption{Weighted triangulations of $\displaystyle\left[\frac{8}{11}\right]_q^{\sharp}$.}
    
    \label{Fig.9-1}
\end{figure}

\begin{figure}[H]
\begin{center}
\begin{tikzpicture}[scale=5]

\draw  (0,0) node[below]{$\displaystyle\left[0\right]_q^{\sharp}$} -- (.5,.6) node[above]{$\displaystyle\left[\infty\right]_q^{\sharp}$}--(.5,0)node[below]{$\displaystyle\left[\frac{1}{1}\right]_q^{\sharp}$}--(0,0);

\draw  (.5,0)--(1,.6)node[above]{$\displaystyle\left[\frac{2}{1}\right]_q^{\sharp}$}--(.5,.6);

\draw  (.5,0)--(1.5,.6)node[above]{$\displaystyle\left[\frac{3}{2}\right]_q^{\sharp}$}--(1,.6);
\draw  (.5,0)--(1,0)node[below]{$\displaystyle\left[\frac{4}{3}\right]_q^{\sharp}$}--(1.5,.6);
\draw  (1,0)--(2,.6)node[above]{$\displaystyle\left[\frac{7}{5}\right]_q^{\sharp}$}--(1.5,.6);
\draw  (1,0)--(2.5,.6)node[above]{$\displaystyle\left[\frac{11}{8}\right]_q^{\sharp}$}--(2,.6);

\draw(.25,0)node[above]{$\displaystyle 1$};
\draw(.5,.3)node[right]{$\displaystyle1$};

\draw(.75,.3)node[right]{$\displaystyle 1$};
\draw(.75,.6)node[above]{$\displaystyle q$};
\draw(1.05,.306)node[above]{$\displaystyle 1$};
\draw(.75,0)node[above]{$\displaystyle 1$};

\draw(1.25,.6)node[above]{$\displaystyle q$};
\draw(1.25,.3)node[left]{$\displaystyle q$};

\draw(1.5,.3)node[above]{$\displaystyle q$};
\draw(1.75,.6)node[above]{$\displaystyle 1$};

\draw(2.25,.6)node[above]{$\displaystyle q^2$};
\draw(1.75,.3)node[below]{$\displaystyle 1$};

    \end{tikzpicture}
    \end{center}
    \caption{Weighted triangulations of $\displaystyle\left[\frac{11}{8}\right]_q^{\sharp}$.}
    
    \label{Fig.10-1}
\end{figure}
}
\end{exam}

\begin{exam}
{\normalfont
The weighted triangulations of $\displaystyle\left[\frac{8}{11}\right]_q^{\flat}=[0,1,2,1,1,1]_q^{\flat}=[[1,4,3]]_q^{\flat}$ and $\displaystyle\left[\frac{11}{8}\right]_q^{\flat}=[1,2,2,1]^{\flat}_q=[[2,2,4]]_q^{\flat}$ can be expressed as Figures \ref{Fig.9} and \ref{Fig.10}. 
\begin{figure}[H]
\begin{center}
\begin{tikzpicture}[scale=5.5]

    \draw  (0,0) node[below]{$\displaystyle\left[0\right]_q^{\flat}$} -- (0,0.6) node[above]{$\displaystyle\left[\infty\right]_q^{\flat}$}--(.5,.6)node[above]{$\displaystyle\left[\frac{1}{1}\right]_q^{\flat}$}--(0,0);

\draw  (0,0)--(.5,0)node[below]{$\displaystyle\left[\frac{1}{2}\right]_q^{\flat}$}--(.5,.6);
\draw  (.5,0)--(1,0)node[below]{$\displaystyle\left[\frac{2}{3}\right]_q^{\flat}$}--(.5,.6);
\draw  (1,0)--(1,.6)node[above]{$\displaystyle\left[\frac{3}{4}\right]_q^{\flat}$}--(.5,.6);
\draw  (1,.6)--(1.5,0)node[below]{$\displaystyle\left[\frac{5}{7}\right]_q^{\flat}$}--(1,0);
\draw  (1.5,0)--(1.5,.6)node[above]{$\displaystyle\left[\frac{8}{11}\right]_q^{\flat}$}--(1,.6);

\draw(.25,0)node[above]{$\displaystyle q^2$};
\draw(.25,.3)node[above]{$\displaystyle q$};
\draw(.25,.6)node[above]{$\displaystyle 1$};
\draw(.5,.3)node[right]{$\displaystyle 1$};
\draw(.75,.0)node[above]{$\displaystyle q$};
\draw(.75,.3)node[above]{$\displaystyle 1$};
\draw(.75,.6)node[above]{$\displaystyle 1$};
\draw(1,.3)node[right]{$\displaystyle q$};
\draw(1.25,.0)node[above]{$\displaystyle q^2$};
\draw(1.25,.3)node[above]{$\displaystyle 1$};
\draw(1.25,.6)node[above]{$\displaystyle 1$};
\draw(1.5,.3)node[right]{$\displaystyle q$};

    \end{tikzpicture}
    \end{center}
    \caption{Weighted triangulations of $\displaystyle\left[\frac{8}{11}\right]_q^{\flat}$.}
    
    \label{Fig.9}
\end{figure}

\begin{figure}[H]
\begin{center}
\begin{tikzpicture}[scale=5]

\draw  (0,0) node[below]{$\displaystyle\left[0\right]_q^{\flat}$} -- (.5,.6) node[above]{$\displaystyle\left[\infty\right]_q^{\flat}$}--(.5,0)node[below]{$\displaystyle\left[\frac{1}{1}\right]_q^{\flat}$}--(0,0);

\draw  (.5,0)--(1,.6)node[above]{$\displaystyle\left[\frac{2}{1}\right]_q^{\flat}$}--(.5,.6);

\draw  (.5,0)--(1.5,.6)node[above]{$\displaystyle\left[\frac{3}{2}\right]_q^{\flat}$}--(1,.6);
\draw  (.5,0)--(1,0)node[below]{$\displaystyle\left[\frac{4}{3}\right]_q^{\flat}$}--(1.5,.6);
\draw  (1,0)--(2,.6)node[above]{$\displaystyle\left[\frac{7}{5}\right]_q^{\flat}$}--(1.5,.6);
\draw  (1,0)--(2.5,.6)node[above]{$\displaystyle\left[\frac{11}{8}\right]_q^{\flat}$}--(2,.6);

\draw(.25,0)node[above]{$\displaystyle q$};
\draw(.5,.3)node[right]{$\displaystyle1$};

\draw(.75,.3)node[right]{$\displaystyle q$};
\draw(.75,.6)node[above]{$\displaystyle 1$};
\draw(1.05,.306)node[above]{$\displaystyle q^2$};
\draw(.75,0)node[above]{$\displaystyle q^3$};

\draw(1.25,.6)node[above]{$\displaystyle 1$};
\draw(1.25,.3)node[left]{$\displaystyle 1$};

\draw(1.5,.3)node[above]{$\displaystyle q$};
\draw(1.75,.6)node[above]{$\displaystyle 1$};

\draw(2.25,.6)node[above]{$\displaystyle 1$};
\draw(1.75,.3)node[below]{$\displaystyle q^2$};

    \end{tikzpicture}
    \end{center}
    \caption{Weighted triangulations of $\displaystyle\left[\frac{11}{8}\right]_q^{\flat}$.}
    
    \label{Fig.10}
\end{figure}
}
\end{exam}



\section{Jones polynomial and left $q$-rational numbers}\label{Sec4}
Following \cite[Proposition 1.2 (b)]{KR}, \cite[Proposition A.1]{MO1} and \cite[Theorem A.3]{BBL}, we can obtain the relationship between left and right $q$-deformed rational numbers and Jones polynomials. In this section, we give a new proof of Theorem A.3 in \cite{BBL} without the homological argument. 

\indent
For a rational number $\displaystyle \alpha=[[c_1,\ldots,c_k]] \in \mathbb{Q}_{>1}$, following \cite{BBL}, we suppose that $\displaystyle V_{\alpha}(q)$ is the
Jones polynomial associated with the rational knot parametrized by $\alpha$, and $\displaystyle \left\vert V_{\alpha}(q) \right\vert$ denote the polynomial obtained by making each coefficient positive. Following \cite{MO1}, let $J_{\alpha}(q)$ be the \textit{normalized Jones polynomial} associated with the rational knot parametrized by $\alpha$. The next lemma can be checked by \cite[Proposition A.1]{MO1} and Theorem \ref{THM-F-S-R}.

\begin{lem}\label{R-F-Jones}
For a rational number $\displaystyle \alpha=[[c_1,\ldots,c_k]]$ which is the Farey sum of $\beta$ and $\gamma$ defined by \eqref{BETA} and \eqref{GAMMA}. 
Then one has
\[
\displaystyle
J_{\alpha}(q)=J_{\beta}(q)+q^{c_k-1}J_{\gamma}(q).
\]
\end{lem}

\indent
Following \cite{BBL}, the sequence of
coefficients of the normalized Jones polynomial $J_{\alpha}(q)$ is just the reverse of the sequence of coefficients of the Jones polynomial $\displaystyle \left\vert V_{\alpha}(q) \right\vert$, then the equation 
\[
\displaystyle \left\vert V_{\alpha}(q) \right\vert=\mathcal{R}^{\flat}_{\alpha}(q)
\]
will be proved by showing the next theorem.

\begin{thm}\label{J-L-Q-R}
For a rational $\displaystyle \alpha=[[c_1,\ldots,c_k]]$, the normalized Jones polynomial of $\alpha$ satisfies the following formula:
\[
\displaystyle
J_{\alpha}(q)=q^m\mathcal{R}^{\flat}_{\alpha}(q^{-1}),
\]
where $\displaystyle m=\mathrm{deg}(\mathcal{R}^{\flat}_{\alpha}(q))=\sum_{j=1}^{k}c_j-k+1$. 
\end{thm}

\medskip

\begin{proof}
For the rational $\displaystyle \alpha=[[c_1,\ldots,c_k]]$, it is easy to check the case of $k=1$ and $k=2$.
We consider the following induction hypothesis:
\begin{equation}\label{H}
\displaystyle
J_{[[c_1,\ldots,c_i]]}(q)=q^m\mathcal{R}^{\flat}_{[[c_1,\ldots,c_i]]}(q^{-1}), \ \ \ \ \text{for} \ 1\leq i \leq k.
\end{equation}
We prove that 
\[
\displaystyle
J_{\alpha^{\prime}}(q)=q^{\mathrm{deg}(\mathcal{R}^{\flat}_{\alpha^{\prime}}(q))}\mathcal{R}^{\flat}_{\alpha^{\prime}}(q^{-1}),
\]
where $\displaystyle\alpha^{\prime}=[[c_1,\ldots,c_k,c_{k+1}]]$.

For the case of $c_k=2$, we have \[\displaystyle \alpha^{\prime}=\beta^{\prime} \# \alpha,\]
where $\displaystyle
\beta=\left\{    
\begin{array}{
    cl}
    \ [[c_1,\ldots,c_{l}-1]] &  \text{for}  \  \   \  \ c_k=c_{k-1}=\cdots=c_{l+1}=2, \ c_{l}>2 ,\ 1\leq l\leq k,   \\ & \\
    \ [[1]]  &  \text{for}  \ \ k=1 , \ c_k=2. \\
\end{array} \right.$
\\
\
\\
Suppose that \[\displaystyle
m_1=\mathrm{deg}(\mathcal{R}^{\flat}_{\beta^{\prime}}(q))=\left\{    
\begin{array}{
    cl}
    \ \displaystyle\sum^{l}_{j=1}c_j-l &  \text{for}  \  \   \  \ c_{k+1}=c_{k}=\cdots=c_{l+1}=2, \ c_{l}>2 ,\ 1\leq l\leq k+1,   \\ & \\
    \ 1  &  \text{for}  \ \  k=1 , \ c_k=2, \\
\end{array} \right.\] 
then by the induction hypothesis \eqref{H}, Lemma \ref{R-F-Jones} and Theorem \ref{THM-F-S-L}, it follows that
\begin{align}
\displaystyle
J_{\alpha^{\prime}}(q)&=J_{\beta^{\prime}}(q)+qJ_{\alpha}(q) \notag \\
&=q^{m_1}\mathcal{R}^{\flat}_{\beta^{\prime}}(q^{-1})+q^{m+1}\mathcal{R}^{\flat}_{\alpha}(q^{-1}) \notag \\
&=q^{m+1}(q^{-(m-m_1+1)}\mathcal{R}^{\flat}_{\beta^{\prime}}+\mathcal{R}^{\flat}_{\alpha}(q^{-1})) \notag \\
&=q^{m+1}\mathcal{R}^{\flat}_{\alpha^{\prime}}(q^{-1})),\notag
\end{align}
where $\displaystyle m+1=\sum^{k}_{j=1}c_j-k+2=\mathrm{deg}(\mathcal{R}^{\flat}_{\alpha^{\prime}}(q)).$

\medskip

Now we assume $\alpha^{\prime}=[[c_1,\ldots,c_k,c]]$ for some $c\geq2$. 
We set the following induction hypothesis:
\begin{equation}\label{H1}
\displaystyle
J_{\alpha^{\prime}}(q)=q^{m^{\prime}}\mathcal{R}^{\flat}_{\alpha^{\prime}}(q^{-1}),
\end{equation}
where $\displaystyle m^{\prime}=\mathrm{deg}(\mathcal{R}^{\flat}_{\alpha^{\prime}}(q))=\sum^{k}_{j=1}c_j-k+c.$

Suppose that $\alpha^{\prime\prime}=[[c_1,\ldots,c_k,c+1]]$, then we have $\displaystyle
\alpha^{\prime\prime}=\alpha^{\prime}\#\alpha.
$
By the induction hypothesis \eqref{H}, \eqref{H1}, Lemma \ref{R-F-Jones} and Theorem \ref{THM-F-S-L}, it follows that
\begin{align}
\displaystyle
J_{\alpha^{\prime\prime}}(q)&=J_{\alpha^{\prime}}(q)+q^cJ_{\alpha}(q) \notag \\
&=q^{m^{\prime}}\mathcal{R}^{\flat}_{\alpha^{\prime}}(q^{-1})+q^{m+c}\mathcal{R}^{\flat}_{\alpha}(q^{-1}) \notag \\
&=q^{m+c}(q^{-1}\mathcal{R}^{\flat}_{\alpha^{\prime}}+\mathcal{R}^{\flat}_{\alpha}(q^{-1})) \notag \\
&=q^{m+c}\mathcal{R}^{\flat}_{\alpha^{\prime\prime}}(q^{-1})),\notag
\end{align}
where $\displaystyle m+c=\sum^{k}_{j=1}c_j-k+1+c=\mathrm{deg}(\mathcal{R}^{\flat}_{\alpha^{\prime\prime}}(q)).$
\end{proof}

\medskip

\begin{exam}
{\normalfont
For $\displaystyle \frac{9}{4}=[[3,2,2,2]]$, since
\[
\displaystyle
\mathrm{deg}(\mathcal{R}^{\flat}_{\frac{9}{4}}(q))=3+2+2+2-4+1=6,
\]
and
\[
\displaystyle 
\left[\frac{9}{4}\right]^{\flat}_q=\frac {1+q+{q}^{2}+2\,{q}^{3}+2\,{q}^{4}+{q}^{5}+{q}^{6}}{1+{q}^{2}+{
q}^{3}+{q}^{4}},
\]
then by Theorem \ref{J-L-Q-R}, we have
\[
\displaystyle
J_{\frac{9}{4}}(q)=q^6\mathcal{R}^{\flat}_{\frac{9}{4}}(q^{-1})=1+q+2q^2+2q^3+q^4+q^5+q^6.
\]
}
\end{exam}

\medskip

\section{Relationship to $2$-Calabi–Yau category associated to the $A_2$ quiver}\label{Sec5}

A relation between $q$-deformed rational numbers and homology algebra is given in \cite{BBL}. In this section, we first briefly recall the relevant definitions and conclusions. We derive a homological interpretation of $q$-Farey sum of $q$-deformed rational numbers by combining Theorems 3.7 and 3.8 in \cite{BBL} with Theorems \ref{THM-F-S-R} and \ref{THM-F-S-L} in Section \ref{Sec3}. In addition, we consider any continued fraction expansion of a real quadratic irrational number of purely cyclic type and give its homological interpretation.

\medskip

A relation between $q$-deformed rational numbers and homology algebra is given in \cite{BBL}. In this chapter, we briefly recall the relevant definitions and conclusions (c.f. \cite{Bri, BBL, BBL2, BDL2, ST01} for more details in this chapter).

\subsection{$2$-Calabi–Yau category associated to the $A_2$ quiver and spherical objects}

\medskip

We fix a field $\Bbbk$ which is algebraic closed and $\mathrm{char}(\Bbbk)=0$.
Consider $DA_2$ the double of $A_2$ quiver as
\[
\begin{tikzcd}
1 \arrow[r,bend left] & 2 \arrow[l,bend left,swap],
\end{tikzcd}
\] where $1$ and $2$ are vertices.
\begin{defn}[Zigzag algebra of $A_2$]
{\normalfont
The zigzag algebra $Z_2$ of $A_2$ is the quotient of the path algebra
$\Bbbk(DA_2)$ by the two-sided ideal generated by all length three paths.
}
\end{defn}
 We regard $Z_2$ as a differential graded algebra by assuming that the grading is the path length and the differential is zero. We denote the homotopy category of differential graded modules over $Z_2$ by $\mathcal{H}_2$, and denote by $\mathcal{D}_2$ the derived category of differential graded modules over $Z_2$ (by inverting quasi-isomorphisms in the $\mathcal{H}_2$). Note that $\mathcal{D}_2$ is not the derived category of complexes of graded $Z_2$-module as an ordinary algebra (c.f. \cite[Section 4]{ST01}). 
For $i=1,2$, we denote the differential graded module $Z_2(i)$ by $P_i$, where $(i)$ is the trivial path. By an ambusing notion, we will denote $P_i$ as an object of $\mathcal{D}_2$.
\begin{defn}
{\normalfont
Let $\mathcal{C}_2$ be the full triangulated subcategory of $\mathcal{D}_2$ which is the extension closure of $\left\{ P_1^{\oplus m_1}[s_1]\bigoplus P_2^{\oplus m_2}[s_2] : m_i \in \mathbb{Z}_{\geq 0}, s_i \in \mathbb{Z}, i=1,2 \right\}.$
}
\end{defn}

\begin{thm}[\cite{Bri, BBL, BBL2, BDL2, ST01}]
For the $\mathcal{C}_2$, we have the following facts.\\
\begin{enumerate}
\item[$(i)$]  $\mathcal{C}_2$ is a finite type linear triangulated category;
\item[$(ii)$] $\mathcal{C}_2$ is a $2$-Calabi–Yau category;
\item[$(iii)$] $P_1$ and $P_2$ are spherical objects in the finite type $\mathcal{C}_2$, that is for $i,j \in {1,2}$ and any integer~$m$, one has
\[
\displaystyle
\mathrm{Hom}(P_i,P_j[m])=\left\{    
\begin{array}{
    cl}
\Bbbk & \qquad \text{for} \ m=0,2 \  \text{and} \ i=j, \\
\Bbbk  & \qquad \text{for} \ m=1 \  \text{and} \ i \neq j, \\
0 & \qquad \text{for}\ otherwise. \\
\end{array} \right.
\]
\end{enumerate}

\end{thm}
There is a unique morphism $\varphi_{12} : P_1[-1] \rightarrow P_2$, and also $\varphi_{21}: P_2[-1]\rightarrow P_1$, and we denote the cones of $\varphi_{12}$ and $\varphi_{21}$ by $P_{12}$ and $P_{21}$, respectively. We note that $P_1$, $P_2$, $P_{12}$ and $P_{21}$ are indecomposable spherical objects of $\mathcal{C}_2$.\\

The extension closure of $P_1$ and $P_2$ is a heart of a bounded $t$-structure of $\mathcal{C}_2$ (c.f. \cite[Section 10.1]{KMSP}, for associated definitions). We call this heart is a standard heart, and denote it by $\heartsuit_{std}$. Note that $P_1$ and $P_2$ are simple objects in $\heartsuit_{std}$, and $\heartsuit_{std}$ is a module category of the preprojective algebra of type $A_2$. We have a exact sequence
\begin{center}
\begin{tikzpicture}
\node (A) at (-1,0) {$0$}; \node (B) at (1,0) {$P_1$};
\node (C) at (3,0) {$P_{21}$};\node (D) at (5,0) {$P_2$};
\node (E) at (7,0) {$0.$};

\draw[->,>=stealth] (A) --(B);
\draw[->,>=stealth] (B) --(C);
\draw[->,>=stealth] (C) --(D);
\draw[->,>=stealth] (D) --(E);

\end{tikzpicture}
\end{center} 

\medskip
\subsection{Spherical twist functors on $\mathcal{C}_2$}

\begin{defn}[Spherical twist functors on $\mathcal{C}_2$ \cite{ST01}]
{\normalfont
Let $X$ is a spherical object in $\mathcal{C}_2$. The spherical twist functor along $X$ is an autoequivalence on $\mathcal{C}_2$ as follows.
\[
\displaystyle
\sigma_{P_i}(X):=\mathrm{Cone}\left(\mathrm{hom}^{\bullet}(P_i,X)\otimes_{\Bbbk}P_i \xrightarrow{\mathrm{ev}} X \right)
\]
where $\mathrm{hom}^{\bullet}(P_i,X)$ is a complex which for $k\in \mathbb{Z}$ the $k$-th term is defined by 
$$
\mathrm{hom}^{k}(P_i,X):=\bigoplus_{j\in \mathbb{Z}}\mathrm{Hom}(P_{i}[j],X[j-k]).
$$
We also have a inverse of $\sigma_{P_i}$ which is defined as
\[
\displaystyle
\sigma^{-1}_{P_i}(X):=\mathrm{Cone}\left(X \xrightarrow{\mathrm{ev}}\mathrm{Lin}^{\bullet}\left(\mathrm{hom}^{\bullet}(P_i,X), P_i \right) \right)\left[-1\right].
\]
}
\end{defn}
Let $\mathbb{S_O}$ be the set of spherical objects of $\mathcal{C}_2$.
\begin{prop}[Braid relation \cite{ST01}]\label{BRGST}
The group of spherical twist functor which generated by $\sigma_{P_1}$, $\sigma_{P_2}$ is isomorphic to $B_3$. Moreover for any $X \in \mathbb{S}$ and for any spherical twist functors $\sigma$ on $\mathcal{C}_2$, one has $\sigma(X) \in \mathbb{S_O}$.  
\end{prop}

Henceforth, without causing confusion, we simply denote $\sigma_{P_1}$ and $\sigma_{P_2}$ as $\sigma_1$ and $\sigma_2$, respectively.\\

For an irredicuible fraction 
$\alpha=[a_1,a_2,\ldots,a_{2m}]$, let the spherical object corresponding to $\alpha$ be
\begin{equation}\label{Xalpha}
X_{\alpha}:=\sigma_1^{-a_1}\sigma_2^{a_2}\cdots \sigma_1^{-a_{2m-1}}\sigma_2^{a_{2m}}P_1.
\end{equation}
By Proposition \ref{BRGST}, $X_{\alpha}$ belongs to $\mathbb{S_O}$.
\subsection{Bridgeland stability conditions}
\begin{defn}[\cite{Bri}]
{\normalfont
A stability condition on a triangulated category $\mathscr{T}$ is specified by two compatible structures $\tau:=(\mathscr{P},Z)$, where $\mathscr{P}$ is a slicing which consists of ablian subcategories $\mathscr{P(\phi)}$ of $\mathscr{T}$ for each real number $\phi$, and $Z$ is called the central charge which is a homomorphism of additive groups from the Grothendieck group $K_0(\mathscr{T})$ to the complex numbers $\mathbb{C}.$\\
The slicing and the central charge satisfy the following conditions:
\begin{enumerate}
\item[(i)]  $\mathscr{P}(\phi+1)=\mathscr{P(\phi)}[1]$;
\item[(ii)] For $\phi>\psi$, if $X \in Ob\mathscr{P(\phi)}$, $Y \in Ob\mathscr{P(\psi)}$, then $\mathrm{Hom}(X,Y)=0$;
\item[(iii)] For any nonzero object $X$ in $\mathscr{T}$, there is a Harder–Narasimhan (HN) filtration of $X$, where each triangle in Figure \ref{Fig.HN-F} is exact and each $Y_i \in \mathscr{P}(\phi_{i})$ for $\phi_1>\phi_2>\cdots>\phi_n$;

\begin{figure}[H]
\begin{center}
\tikzset{cross/.style={preaction={-,draw=white,line width=6pt}}}
\begin{tikzpicture}
\node (0=) at (-1.8,2) {$0=$};
\node (X_0) at (-1,2) {$X_0$}; \node (X_1) at (1,2) {$X_1$};
\node (X_2) at (3,2) {$X_2$}; \node (d) at (5,2) {$\cdots$};
\node (X_{n-2}) at (7,2) {$X_{n-2}$};
\node (X_{n-1}) at (9,2) {$X_{n-1}$};
\node (X_{n}) at (11,2) {$X_{n}$};
\node (=X) at (11.8,2) {$=X$};

\node (Y_n) at (10,0) {$Y_n$};
\node (Y_{n-1}) at (8,0) {$Y_{n-1}$};
\node (Y_2) at (2,0) {$Y_2$};
\node (Y_1) at (0,0) {$Y_1$};

\draw[->,>=stealth] (X_0) --(X_1);
\draw[->,>=stealth] (X_1) --(X_2);
\draw[->,>=stealth] (X_2) --(d);
\draw[->,>=stealth] (d) --(X_{n-2});
\draw[->,>=stealth] (X_{n-2}) --(X_{n-1});
\draw[->,>=stealth] (X_{n-1}) --(X_{n});

\draw[->,>=stealth] (X_{n}) --(Y_n);
\draw[dashed,->,>=stealth] (Y_n) --(X_{n-1});

\draw[->,>=stealth] (X_{n-1}) --(Y_{n-1});
\draw[dashed,->,>=stealth] (Y_{n-1}) --(X_{n-2});

\draw[->,>=stealth] (X_2) --(Y_2);
\draw[dashed,->,>=stealth] (Y_2) --(X_1);

\draw[->,>=stealth] (X_1) --(Y_1);
\draw[dashed,->,>=stealth] (Y_1) --(X_0);

\end{tikzpicture}

\end{center}
\caption{The Harder–Narasimhan filtrations of X}
    
    \label{Fig.HN-F}
\end{figure}

\item[(iv)] For any $X \in Ob\mathscr{P(\phi)}$, there is some positive real number $m$ such that $$Z([X]) = me^{i\pi\phi}.$$
\end{enumerate}
}
\end{defn}
\begin{defn}[\cite{BBL}]
{\normalfont
A stability condition on $\mathcal{C}_2$ is standard if there is a positive real number $\phi$, such that $\mathscr{P}([\phi,\phi +1))=\heartsuit_{std}$. More specifically, $\mathscr{P}([\phi,\phi +1))$ is formed by $P_1$ and $P_2$ on the spot. 
}
\end{defn}

A standard stability condition $\tau$ gives rise to a special filtration of $X_{\alpha}$ (see Figure \ref{Fig.HN-F_X_alpha}),

\begin{figure}[H]
\begin{center}
\tikzset{cross/.style={preaction={-,draw=white,line width=6pt}}}
\begin{tikzpicture}
\node (0=) at (-2,2) {$0=$};
\node (X_0) at (-1,2) {$(X_{\alpha})_0$}; \node (X_1) at (1,2) {$(X_{\alpha})_1$};
\node (X_2) at (3,2) {$(X_{\alpha})_2$}; \node (d) at (5,2) {$\cdots$};
\node (X_{n-2}) at (7,2) {$(X_{\alpha})_{n-2}$};
\node (X_{n-1}) at (9.2,2) {$(X_{\alpha})_{n-1}$};
\node (X_{n}) at (11.2,2) {$(X_{\alpha})_{n}$};
\node (=X) at (12.4,2) {$=X_{\alpha}$};

\node (Y_n) at (10,0) {$Y_n$};
\node (Y_{n-1}) at (8,0) {$Y_{n-1}$};
\node (Y_2) at (2,0) {$Y_2$};
\node (Y_1) at (0,0) {$Y_1$};

\draw[->,>=stealth] (X_0) --(X_1);
\draw[->,>=stealth] (X_1) --(X_2);
\draw[->,>=stealth] (X_2) --(d);
\draw[->,>=stealth] (d) --(X_{n-2});
\draw[->,>=stealth] (X_{n-2}) --(X_{n-1});
\draw[->,>=stealth] (X_{n-1}) --(X_{n});

\draw[->,>=stealth] (X_{n}) --(Y_n);
\draw[dashed,->,>=stealth] (Y_n) --(X_{n-1});

\draw[->,>=stealth] (X_{n-1}) --(Y_{n-1});
\draw[dashed,->,>=stealth] (Y_{n-1}) --(X_{n-2});

\draw[->,>=stealth] (X_2) --(Y_2);
\draw[dashed,->,>=stealth] (Y_2) --(X_1);

\draw[->,>=stealth] (X_1) --(Y_1);
\draw[dashed,->,>=stealth] (Y_1) --(X_0);

\end{tikzpicture}

\end{center}
\caption{The Harder–Narasimhan filtrations of $X_{\alpha}$}
    
    \label{Fig.HN-F_X_alpha}
\end{figure}
\noindent
where each subquotient (called $\tau$-HN-factor) $Y_j$ is isomorphic to $P_{\nu}[i]$ ($j=0,\ldots,n$, $i \in \mathbb{Z}$, $\nu \in \{1,2,12,21\}$).



\subsection{Two kinds of functional $\mathrm{occ}_q$ and $\overline{\mathrm{hom}}_q$}

Fix a standard stability condition $\tau$, let $[P_i,P_j]$ denote the set of all objects of $\mathcal{C}_2$ whose $\tau$-Harder-Narasimhan filtration factors are shifts of either $P_i$ or $P_j$ where $i,j=1,2,12,21$ and $i\neq j$. By \cite[Proposition 4.3]{BBL}, each spherical object $X \in \mathbb{S_O}$ must belong to one of $[P_2,P_{21}]$, $[P_{21},P_1]$, $[P_1,P_{12}]$, $[P_{12},P_2]$. 
Since we only consider the case $\alpha\in\mathbb{Q} \cap (0,\infty)$, then by \cite[Figure 6 in Section 4.3]{BBL}, it must have $X \in [P_2,P_{21}]\cup[P_{21},P_1]$. 

\begin{defn}[Counting functional]
{\normalfont
The counting functional $\pi_{\nu}(X)$ of $\tau$-HN-factors is defined by an element of the Laurent polynomial ring $\mathbb{Z}[q,q^{-1}]$ as follows.
$$
\displaystyle
\pi_{\nu}(X):=\sum_{i\in \mathbb{Z}}\eta_{\nu}q^i,
$$
where $\nu \in \{1,2,12,21\}$ and $\eta_{\nu}$ is the number of $\tau$-HN-factors which are isomorphic to $P_{\nu}$. 
}
\end{defn}

\begin{defn}[{\cite[$\mathrm{occ}_q$ and $\overline{\mathrm{hom}}_q$]{BBL}}]
{\normalfont
For any $X,Y \in \mathbb{S_O}$,  there are two kinds of functionals $\mathbb{S_O}\times\mathbb{S_O}\rightarrow~\mathbb{Z}[q^{\pm}]$ denoted by $\mathrm{occ}_q$ and $\overline{\mathrm{hom}}_q$, respectively, which are defined as follows.
\[
\displaystyle
\mathrm{occ}_q(P_1,X):=\pi_2(X)+\pi_{12}(X)+\pi_{21}(X),
\]
\[
\displaystyle
\mathrm{occ}_q(P_2,X):=\pi_1(X)+\pi_{12}(X)+\pi_{21}(X),
\]
\[
\displaystyle
\overline{\mathrm{hom}}_q(X,Y):=\left\{    
\begin{array}{c@{\hspace{0.5cm}}l}
  q^{n}(q^{-2}-q^{-1})  \  &  \text{if} \ \  Y \cong X[n], \\ & \\
    \ \sum_{n\in \mathbb{Z}}\mathrm{dim}_{\Bbbk}\mathrm{Hom}(X,Y[n])q^{-n}  &     \text{otherwise}. \\
\end{array} \right.
\]
}
\end{defn}

\medskip

Then we have the next theorem which gives a relationship between the $q$-deformed rational numbers and homological algebra.

\begin{thm}[Bapat, Becker and Licata {\cite[A part of Theorems 4.7 and 4.8]{BBL}}]\label{3.7_and_3.8}

Consider a non-negative rational number $\alpha=[a_1,a_2,\ldots,a_{2m}]$. Suppose that $$
X_{\alpha}=\sigma_1^{-a_1}\sigma_2^{a_2}\cdots \sigma_1^{-a_{2m-1}}\sigma_2^{a_{2m}}P_1,
$$ then we have
\[
\displaystyle 
[\alpha]_q^{\sharp}=\frac{\mathrm{occ}_q(P_2,X_\alpha)}{\mathrm{occ}_q(P_1,X_\alpha)},
\]
and
\[
\displaystyle
[\alpha]_q^{\flat}=\frac{\overline{\mathrm{hom}}_q(X_\alpha,P_2)}{q\overline{\mathrm{hom}}_q(X_\alpha,P_1)}.
\]
\end{thm}

 
\medskip

\subsection{A homological interpretation of the $q$-deformed Farey sum}
Consider a positive rational number $\alpha=[a_1,a_2,\ldots,a_{2m}]=[[c_1,\ldots,c_k]]$. Since Theorems \ref{THM-F-S-R} and \ref{THM-F-S-L} are based on $[[c_1,\ldots,c_k]]$, we first make a simple formal transformation of \eqref{Xalpha}. By direct computation, $S$ in Section \ref{Sec2} is represented by $S=\sigma_1\sigma_2\sigma_1$ in $\mathrm{PSL}_2(\mathbb Z)$, and then we have
\begin{align}
\displaystyle
\sigma_1^{-a_1}\sigma_2^{a_2}\cdots \sigma_1^{-a_{2m-1}}\sigma_2^{a_{2m}}&=\sigma_1^{-c_1}S\sigma_1^{-c_2}S\cdots \sigma_1^{-c_k}S\sigma_1^{-1} \notag \\
 &= \sigma_1^{-c_1+1}\sigma_2\sigma_1^{-c_2+2}\sigma_2\sigma_1^{-c_3+2}\sigma_2\cdots\sigma_1^{-c_k+2}\sigma_2, \notag
\end{align}
in $\mathrm{PSL}_2(\mathbb Z)$. Thus, we also have
$$
\sigma_1^{-a_1}\sigma_2^{a_2}\cdots \sigma_1^{-a_{2m-1}}\sigma_2^{a_{2m}}=\sigma_1^{-c_1+1}\sigma_2\sigma_1^{-c_2+2}\sigma_2\sigma_1^{-c_3+2}\sigma_2\cdots\sigma_1^{-c_k+2}\sigma_2,
$$
in $B_3$. Hence, we have 
\begin{equation}\label{XXalpha}
\displaystyle
X_{\alpha}=\sigma_1^{-c_1+1}\sigma_2\sigma_1^{-c_2+2}\sigma_2\sigma_1^{-c_3+2}\sigma_2\cdots\sigma_1^{-c_k+2}\sigma_2P_1.
\end{equation}
If $\alpha$ is the Farey sum of $\beta$ and $\gamma$ defined by \eqref{BETA} and \eqref{GAMMA}, then by appling the above transformation to $\beta$ and $\gamma$, it follows that 
\begin{equation}\label{XXbeta}
\displaystyle
X_{\beta}=\left\{    
\begin{array}{
    c@{}l}
  \sigma_1^{-c_1+1}\sigma_2\sigma_1^{-c_2+2}\sigma_2\sigma_1^{-c_3+2}\sigma_2\cdots & \sigma_1^{-c_{l-1}+2}\sigma_2\sigma_1^{-c_l+3}\sigma_2P_1   \   \\ 
  & \text{for} \ \ \    c_k=c_{k-1}=\cdots=c_{l+1}=2, \ c_{l}>2 ,\ 1\leq l\leq k,   \\ & \\
    \ \sigma_2P_1  &  \text{for}  \ \ k=1 , \ c_k=2,  \\ & \\
     \ P_2  &  \text{for}  \ \ c_k=c_{k-1}=\cdots=c_2=2 , \ c_1=1,  \\
\end{array} \right.
\end{equation}
and
\begin{equation}\label{XXgamma}
\displaystyle
X_{\gamma}=\left\{    
\begin{array}{
    cl}
  \sigma_1^{-c_1+1}\sigma_2\sigma_1^{-c_2+2}\sigma_2\sigma_1^{-c_3+2}\sigma_2\cdots\sigma_1^{-c_{k-1}+2}\sigma_2P_1  \  &  \text{for}  \  \   k\geq 2, \\ & \\
    \ P_1  &  \text{for}  \ \ k=1. \\
\end{array} \right.
\end{equation}

Note that $\mathrm{occ}_q(P_2,X_{\alpha})$ and $\mathrm{occ}_q(P_1,X_{\alpha})$ are not the numerator and denominator of the right $q$-deformed rational number, respectively. In fact, they diﬀer by some power of $q$ (c.f. \cite[The proof of Theorems 4.7 and 4.8]{BBL}). Hence, we first determine the power of $q$ as follows.

\begin{thm}\label{key1}
For the right $q$-deformed rational number $\displaystyle\left[\alpha\right]^{\sharp}_q=\frac{\mathcal{R}^{\sharp}_{\alpha}(q)}{\mathcal{S}^{\sharp}_{\alpha}(q)}$, we have
\[
\displaystyle
\mathcal{R}^{\sharp}_{\alpha}(q)=q^{k-1}\mathrm{occ}_q(P_2,X_{\alpha}),\qquad \mathcal{S}^{\sharp}_{\alpha}(q)=q^{k-1}\mathrm{occ}_q(P_1,X_{\alpha}).
\]
\end{thm}

\begin{proof}
For a polynomial $f(q) \in \mathbb Z [q^{-1},q] $, let $\displaystyle md(f(q)):=\mathrm{min}\left\{ i : f(q)=\sum_{i}\rho_i q^i, \rho_i\neq 0 \right\}.$ One has  $md(\mathcal{R}^{\sharp}_{\alpha}(q))=0.$
We consider the $\tau$-HN-factor $Y_j$ of $X_{\alpha}$, and let $$\iota(X_{\alpha}):=\mathrm{min} \left\{ i : P_{\nu}[i]\cong Y_j, \nu=1,2,12,21 \right\}.$$ Then by the definition of $\mathrm{occ}_q$, we have $md(\mathrm{occ}_q(P_2,X_{\alpha}))=md(\mathrm{occ}_q(P_1,X_{\alpha}))=\iota(X_{\alpha})$. It can be checked that for $n_1,n_2,l_1,l_2 \in \mathbb Z_{\geq 0},$
$$\qquad \ \iota(P_1)=0, \quad \iota(\sigma_2P_1)=0, \quad \iota(\sigma_2^{n_1}P_1)=n_1-1,\quad  \iota(\sigma_1^{l_1}\sigma_2^{n_1}P_1)=n_1-1, $$
$$\iota(\sigma_2^{n_2}\sigma_1^{l_1}\sigma_2^{n_1}P_1)=n_1+n_2-1,\quad\iota(\sigma_1^{l_2}\sigma_2^{n_2}\sigma_1^{l_1}\sigma_2^{n_1}P_1)=n_1+n_2-1,$$
by direct computation. For \eqref{XXalpha}, repete the above computation, we have $\iota(X_{\alpha})=k-1.$ Hence 
we have  
$$
\mathcal{R}^{\sharp}_{\alpha}(q)=q^{k-1}\mathrm{occ}_q(P_2,X_{\alpha}),\qquad \mathcal{S}^{\sharp}_{\alpha}(q)=q^{k-1}\mathrm{occ}_q(P_1,X_{\alpha}).
$$
\end{proof}

Similarly, we note that $\overline{\mathrm{hom}}_q(X_{\alpha},P_2)$ and $\overline{\mathrm{hom}}_q(X_{\alpha},P_1)$ are not the numerator and denominator of the right $q$-deformed rational number, respectively. However, we have

\begin{thm}\label{key2}
For the left $q$-deformed rational number $\displaystyle\left[\alpha\right]^{\flat}_q=\frac{\mathcal{R}^{\flat}_{\alpha}(q)}{\mathcal{S}^{\flat}_{\alpha}(q)}$, we have
\[
\displaystyle
\mathcal{R}^{\flat}_{\alpha}(q)=q^{\sum^{k}_{j=1}(c_j-2)+1}\overline{\mathrm{hom}}_q(X_{\alpha},P_2),\qquad \mathcal{S}^{\flat}_{\alpha}(q)=q^{\sum^{k}_{j=1}(c_j-2)+2}\overline{\mathrm{hom}}_q(X_{\alpha},P_1).
\]
\end{thm}

\begin{proof}
It can be proved by a similar argument of Theorem \ref{key1}, and consider $$\iota^{\prime}(X_{\alpha}):=\mathrm{max} \left\{ i : P_{\nu}[i]\cong Y_j, \nu=1,2,12,21 \right\}$$ and the definition of $\overline{\mathrm{hom}}_q$.
\end{proof}

For a rational number $\displaystyle \alpha=[[c_1,\ldots,c_k]]$ which is the Farey sum of $\beta$ and $\gamma$ defined by \eqref{BETA} and \eqref{GAMMA}, if we consider the spherical objects corresponding $\alpha$, $\beta$, $\gamma$, then by Theorems \ref{THM-F-S-R} and \ref{THM-F-S-L}, we have the following facts. 

\begin{prop}\label{CTFS-RandL}
Considering the spherical objects $X_{\alpha}$, $X_{\beta}$ and $X_{\gamma}$ in $\mathcal{C}_2$, for $i=1,2$, we have the following formula. \\

\begin{enumerate}
\item[$(i)$] If $\alpha \in  \mathbb{Q}_{>1}\setminus \mathbb{Z}_{>1} $ , then
\begin{align}
\displaystyle
\mathrm{occ}_q(P_i,X_{\alpha})&=q^{l-k}\mathrm{occ}_q(P_i,X_{\beta})+q^{c_k-2}\mathrm{occ}_q(P_i,X_{\gamma}),\notag \\
\displaystyle
\overline{\mathrm{hom}}_q(X_{\alpha},P_i)&=q^{2(k-l)-\sum^{k}_{j=l+1}c_j}\overline{\mathrm{hom}}_q(X_{\beta},P_i)+q^{-c_k+2}\overline{\mathrm{hom}}_q(X_{\gamma},P_i); \notag 
\end{align}
\item[$(ii)$] 
If $\alpha \in  \mathbb{Z}_{>1} $ $($i.e. $l=k=1$, and $X_{\gamma}=P_1$ $)$, then
\begin{align}
\displaystyle
\mathrm{occ}_q(P_i,X_{\alpha})&=\mathrm{occ}_q(P_i,X_{\beta})+q^{c_1-1}\mathrm{occ}_q(P_i,P_1),\notag \\
\displaystyle
\overline{\mathrm{hom}}_q(X_{\alpha},P_i)&=\overline{\mathrm{hom}}_q(X_{\beta},P_i)+q^{-\left(\sum^{k}_{j=1}\left(c_j-2\right)+1\right)}\overline{\mathrm{hom}}_q(P_1,P_i); \notag 
\end{align}
\item[$(iii)$] 
If $\alpha \in\left(\mathbb{Q}\cap (0,1)\right) \setminus \mathbb{Z}^{-1}_{>1}$ $($where $\mathbb{Z}^{-1}_{>1}:=\left\{\frac{1}{n} : n \in \mathbb{Z}_{>1}  \right\}$ $)$, then
\begin{align}
\displaystyle
\mathrm{occ}_q(P_i,X_{\alpha})&=q^{l-k-1}\mathrm{occ}_q(P_i,X_{\beta})+q^{c_k-2}\mathrm{occ}_q(P_i,X_{\gamma}),\notag \\
\displaystyle
\overline{\mathrm{hom}}_q(X_{\alpha},P_i)&=q^{2(k-l+1)-\sum^{k}_{j=l+1}c_j}\overline{\mathrm{hom}}_q(X_{\beta},P_i)+q^{-c_k+2}\overline{\mathrm{hom}}_q(X_{\gamma},P_i); \notag 
\end{align}
\item[$(iv)$] 
If $\alpha \in  \mathbb{Z}^{-1}_{>1} $ $($i.e. $X_{\beta}=P_2$ and $c_k-2=0$$)$, then
\begin{align}
\displaystyle
\mathrm{occ}_q(P_i,X_{\alpha})&=q^{-(k-1)}\mathrm{occ}_q(P_i,P_2)+\mathrm{occ}_q(P_i,X_{\gamma}),\notag \\
\displaystyle
\overline{\mathrm{hom}}_q(X_{\alpha},P_i)&=q^{-\sum^{k}_{j=l+1}c_j+3k}\overline{\mathrm{hom}}_q(P_2,P_i)+\overline{\mathrm{hom}}_q(X_{\gamma},P_i). \notag 
\end{align}

\end{enumerate}

\end{prop}

By Corollary \ref{CTFS-RandL},  we can botain a formual Farey sum of the sphercial objects in $\mathcal{C}_2$.

\begin{prop}\label{CTFS-OB}
Considering the spherical objects $X_{\alpha}$, $X_{\beta}$ and $X_{\gamma}$ in $\mathcal{C}_2$,  for $i=1,2$, we have the following formula. \\
\begin{enumerate}
\item[$(i)$] If $\alpha \in  \mathbb{Q}_{>1}\setminus \mathbb{Z}_{>1} $ , then
$$
X_{\alpha}=X_{\beta}[l-k]\oplus X_{\gamma}[c_k-2];
$$
\item[$(ii)$] 
If $\alpha \in  \mathbb{Z}_{>1} $ $($i.e. $l=k=1$, and $X_{\gamma}=P_1$ $)$ , then
$$
X_{\alpha}=X_{\beta}\oplus P_1[c_1-1];
$$
\item[$(iii)$] 
If $\alpha \in\left(\mathbb{Q}\cap (0,1)\right) \setminus \mathbb{Z}^{-1}_{>1}$ $($where $\mathbb{Z}^{-1}_{>1}:=\left\{\frac{1}{n} : n \in \mathbb{Z}_{>1}  \right\}$ $)$, then
$$
X_{\alpha}=X_{\beta}[l-k-1]\oplus X_{\gamma}[c_k-2];
$$
\item[$(iv)$] 
If $\alpha \in  \mathbb{Z}^{-1}_{>1} $ $($i.e. $X_{\beta}=P_2$ and $c_k-2=0$$)$, then
$$
X_{\alpha}=P_2[-(k-1)]\oplus X_{\gamma}.
$$
\end{enumerate}

\end{prop}

\bigskip

\subsection{Real quadratic irrationals with periodic type}
Let $x>1$ be a real quadratic irrational number. Since $[x]_q$ can be written as $\displaystyle [x]_q=\frac{\mathcal{R}+\sqrt{\mathcal{P}}}{\mathcal{S}}$ where $\mathcal{R}$, $\mathcal{P}$, $\mathcal{S}$ $\in \mathbb{Z}[q]$, then we have the following conclusions which give a homological interpretation of $[x]_q$.

\begin{thm}\label{RQIN-occ_q_and_hom_q}
Let $\displaystyle x=[[c_1,\ldots,c_k,c_1,\ldots,c_k,\ldots]]>1$ be a real quadratic irrational number which the continued fraction expansion is purely periodic type. Suppose that $\alpha=[[c_1,\ldots,c_k]]$, and $\gamma=[[c_1,\ldots,c_{k-1}]].$ Considering two spherical objects $X_{\alpha}$ and $X_{\gamma}$ in $\in \mathcal{C}_2$ which are given by \eqref{XXalpha} and \eqref{XXgamma}. Then,  we have
\begin{enumerate}
\item[$(1)$]
\begin{equation}\label{RQIN_occ_q}
\displaystyle
\left[x\right]_q=\frac{\mathcal {A}_1+\mathcal{A}_2 + \sqrt{(\mathcal A_1-\mathcal A_2)^2-4q^{\sum^{k}_{i=1}c_i-3k+2}}}{\mathcal B}, 
\end{equation}
\[
\displaystyle
\mathcal {A}_1=\mathrm{occ}_q(P_2,X_{\alpha}), \
\displaystyle
\mathcal A_2 = q^{c_k-2}\mathrm{occ}_q(P_1,X_{\gamma}), \
\mathcal B = 2\mathrm{occ}_q(P_1,X_{\alpha}),
\]
and
\begin{equation}\label{RQIN_hom_q}
\displaystyle
\left[x\right]_q=\frac{\mathcal{A}^{\prime}_1+\mathcal{A}^{\prime}_2+\sqrt{(\mathcal{A}^{\prime}_1-\mathcal{A}^{\prime}_2)^2-4\text{c}(q)}}{\mathcal{B}^{\prime}},
\end{equation}
where
\begin{align}
\displaystyle
\mathcal{A}^{\prime}_1&=(q-1)\overline{\mathrm{hom}}_q(P_1,X_{\alpha})+q\overline{\mathrm{hom}}_q(P_2,X_{\alpha}),\notag \\
\mathcal{A}^{\prime}_2&=q^{c_k-1}(\overline{\mathrm{hom}}_q(P_1,X_{\gamma})+(1-q)\overline{\mathrm{hom}}_q(P_2,X_{\gamma})), \notag \\
\mathcal{B}^{\prime}&=2q(\overline{\mathrm{hom}}_q(P_1,X_{\alpha})+(1-q)\overline{\mathrm{hom}}_q(P_2,X_{\alpha})),\notag
\end{align}
\[
\displaystyle
\text{c}(q)=q^{\sum^{k}_{i=1}c_i-3k+4}-2q^{\sum^{k}_{i=1}c_i-3k+3}+3q^{\sum^{k}_{i=1}c_i-3k+2}-2q^{\sum^{k}_{i=1}c_i-3k+1}+q^{\sum^{k}_{i=1}c_i-3k};
\]

\item[$(2)$] In particular,  if $x=\displaystyle \frac{c_1+\sqrt{c^2_1-4}}{2}$ $($ $c\geq 3$ $)$,  then
\begin{equation}\label{RQIN_occ_q-2}
\displaystyle
\left[x\right]_q=\frac{\mathrm{occ}_q(P_2,X_{\alpha}),+ \sqrt{\left(\mathrm{occ}_q(P_2,X_{\alpha}\right)^2-4q^{c_1-1}}}{2\mathrm{occ}_q(P_1,X_{\alpha})}.
\end{equation}
\end{enumerate}
\end{thm}

\begin{proof}
We only prove ($1$), and the ($2$) can be prove by the same argument.
By Proposition \ref{N_and_D_and_q-Euler}, we have
\[
\displaystyle
\mathcal{R}^{\sharp}_{\alpha}(q)=E^{\sharp}_k(c_1,\ldots,c_k)_q, \qquad \mathcal{S}^{\sharp}_{\alpha}(q)=E^{\sharp}_{k-1}(c_2,\ldots,c_k)_q; 
\]
\[
\displaystyle
\mathcal{R}^{\sharp}_{\gamma}(q)=E^{\sharp}_{k-1}(c_1,\ldots,c_{k-1})_q, \qquad \mathcal{S}^{\sharp}_{\gamma}(q)=E^{\sharp}_{k-2}(c_2,\ldots,c_{k-1})_q.
\]
On the other hand, we have
\[
\displaystyle
\mathcal{R}^{\sharp}_{\alpha}(q)=q^{k-1}\text{occ}_q(P_2,X_{\alpha}), \qquad \mathcal{S}^{\sharp}_{\alpha}(q)=q^{k-1}\text{occ}_q(P_1,X_{\alpha}); 
\]
\[
\displaystyle
\mathcal{R}^{\sharp}_{\gamma}(q)=q^{k-2}\text{occ}_q(P_2,X_{\gamma}), \qquad \mathcal{S}^{\sharp}_{\gamma}(q)=q^{k-2}\text{occ}(P_1,X_{\gamma}).
\]
By \cite[Proposition 4.3]{LM}, since $[x]_q$ can be written as $\displaystyle [x]_q=\frac{\mathcal{R}+\sqrt{\mathcal{P}}}{\mathcal{S}}$, with
\begin{align}
\displaystyle
\mathcal{R}&=E^{\sharp}_k(c_1,\ldots,c_k)_q+q^{c_k-1}E^{\sharp}_{k-2}(c_2,\ldots,c_{k-1})_q ,\notag \\
\mathcal{P}&=(E^{\sharp}_k(c_1,\ldots,c_k)_q-q^{c_k-1}E^{\sharp}_{k-2}(c_2,\ldots,c_{k-1})_q)^2-4q^{\sum^{k}_{i=1}(c_i-1)} ,\notag \\
\mathcal{S}&=2E^{\sharp}_{k-1}(c_2,\ldots,c_k)_q,\notag
\end{align}
then, by a simple substitution, \eqref{RQIN_occ_q} is proved.\\

On the other hand, by \cite[Lemma 4.13]{BBL}, we can know that the relationship between the  $\mathrm{occ}_q$ and $\overline{\mathrm{hom}}_q$ as follows. 
\[
\displaystyle
\overline{\text{hom}}_q(P_1,X)=q^{-1}\mathrm{occ}_q(P_1,X)+(1-q^{-1})\mathrm{occ}_q(P_2,X),
\]
\[
\displaystyle
\overline{\mathrm{hom}}_q(P_2,X)=(q^{-2}-q^{-1})\mathrm{occ}_q(P_1,X)+q^{-1}\mathrm{occ}_q(P_2,X).
\]
By solving the above two equations on $\mathrm{occ}_q$,  we can obtain the following two equations.
\[
\displaystyle
\mathrm{occ}_q(P_1,X)=\frac{q\overline{\mathrm{hom}}_q(P_1,X)}{q+q^{-1}-1}+\frac{(q-q^2)\overline{\mathrm{hom}}_q(P_2,X)}{q+q^{-1}-1},
\]
\[
\displaystyle
\mathrm{occ}_q(P_2,X)=\frac{(q-1)\overline{\mathrm{hom}}_q(P_1,X)}{q+q^{-1}-1}+\frac{q\overline{\mathrm{hom}}_q(P_2,X)}{q+q^{-1}-1}.
\]
Finally, we substitute these two equations into \eqref{RQIN_occ_q} to obtain \eqref{RQIN_hom_q}.
\end{proof}

\begin{rem}
{\normalfont
For (2) of the above Theorem,  if we take $c_1=3$,  then we have $$X_{\alpha}=P_{21}\oplus P_1[1]\oplus P_1[2]$$
hence, we get the case of the golden number as follows.
\begin{equation}\label{RQIN_occ_q-3}
\displaystyle
\left[\frac{1+\sqrt{5}}{2}\right]_q=\frac{q^{-1}\left(\mathrm{occ}_q(P_2,X_{\alpha}),+ \sqrt{\left(\mathrm{occ}_q(P_2,X_{\alpha}\right)^2-4q^{c_1-1}}-2\mathrm{occ}_q(P_1,X_{\alpha})\right)}{2\mathrm{occ}_q(P_1,X_{\alpha})}.
\end{equation}
}
\end{rem}

\section*{Acknowledgments}
I am greatly indebted to Professor Asilata Bapat for many useful discussions. I also grateful to Professor Akishi Ikeda for providing a notes about Stability conditions on triangulated categories, which helped me learn some basic concepts about it. I am greatly indebted to Professor Michihisa Wakui for many useful disscussions and for the guidance. I thank Professors Takeyoshi Kogiso and Kengo Miyamoto for giving me the opportunity to start this work.

\medskip


\bigskip \bigskip 
Xin Ren: Department of Mathematics, Graduate School of Science, Osaka University, Toyonaka Osaka, 560-0043, Japan. 
\par
\textit{E-mail address}: (1) {\bf ren.xin.sci@osaka-u.ac.jp}; (2) {\bf xinren1213@gmail.com}.

\end{document}